\RequirePackage{fix-cm}
\documentclass[smallextended]{svjour3}       
\smartqed  
\usepackage{braket,amsfonts,hyperref}

\usepackage{array}

\usepackage[caption=false]{subfig}

\usepackage{algorithmic}

\usepackage{graphicx,epstopdf}

\usepackage{amsopn}

\usepackage{xspace}
\usepackage{bold-extra}
\usepackage[most]{tcolorbox}
\usepackage{multirow}

\colorlet{texcscolor}{blue!50!black}
\colorlet{texemcolor}{red!70!black}
\colorlet{texpreamble}{red!70!black}
\colorlet{codebackground}{black!25!white!25}

\lstdefinestyle{siamlatex}{%
  style=tcblatex,
  texcsstyle=*\color{texcscolor},
  texcsstyle=[2]\color{texemcolor},
  keywordstyle=[2]\color{texemcolor},
  moretexcs={cref,Cref,maketitle,mathcal,text,headers,email,url},
}

\tcbset{%
  colframe=black!75!white!75,
  coltitle=white,
  colback=codebackground, 
  colbacklower=white, 
  fonttitle=\bfseries,
  arc=0pt,outer arc=0pt,
  top=1pt,bottom=1pt,left=1mm,right=1mm,middle=1mm,boxsep=1mm,
  leftrule=0.3mm,rightrule=0.3mm,toprule=0.3mm,bottomrule=0.3mm,
  listing options={style=siamlatex}
}

\newtcbinputlisting[use counter=example]{\examplefile}[3][]{%
  title={Example~\thetcbcounter: #2},listing file={#3},#1}

\DeclareTotalTCBox{\code}{ v O{} }
{ 
  fontupper=\ttfamily\color{black},
  nobeforeafter,
  tcbox raise base,
  colback=codebackground,colframe=white,
  top=0pt,bottom=0pt,left=0mm,right=0mm,
  leftrule=0pt,rightrule=0pt,toprule=0mm,bottomrule=0mm,
  boxsep=0.5mm,
  #2}{#1}

\usepackage{booktabs}
\usepackage{color}
\usepackage{amsmath}
\usepackage{amssymb}
\allowdisplaybreaks
\usepackage{algorithm}

\newtheorem{assumption}{Assumption}

\usepackage{endnotes}
\let\footnote=\endnote

\usepackage{bm}

\newcommand*{\xbf}{{\mathbf x}}
\newcommand*{\x}{{\mathbf x}}
\newcommand*{\ubf}{{\mathbf u}}
\newcommand*{\ybf}{{\mathbf y}}
\newcommand*{\y}{{\mathbf y}}
\newcommand*{\sgn}{{\text{sgn}}}
\newcommand*{\zbf}{{\mathbf z}}

\newcommand*{\vbf}{{\mathbf v}}
\newcommand*{\wbf}{{\mathbf w}}

\newcommand*{\z}{{\mathbf z}}

\newcommand*{\R}{{\mathbb R}}
\newcommand*{\bR}{{\mathbb R}}

\newcommand*{\E}{{\mathbb E}}
\newcommand*{\bE}{{\mathbb E}}

\usepackage{ulem}
\newcommand*{\lf}{{\underline f}}
\def\be{\begin{enumerate}}
\def\ee{\end{enumerate}}
\newcommand*{\eh}{{\hat{\epsilon}}}
\newcommand{\sigmainf}{\sigma_{\infty}}
\newcommand{\hrho}{\hat{\rho}}
\usepackage{dsfont}

\def\sO{\mathcal{O}}
\def\hx{\mathbf{\hat x}}
\def\prox{\rm prox}
\def\sN{\mathcal{N}}
\def\hrho{\hat \rho}
\def\conv{\rm conv}
\def\prox{\rm prox}
\newcommand*{\rev}{\color{black}{}}
\def\tf{\tilde{f}}
\def\tg{\tilde{g}}
\newcommand*{\revv}{\color{black}{}}
\newcommand*{\revi}{\color{black}{}}
\def\sH{\mathcal{H}}
\def\sA{\mathcal{A}}
\def\sP{\mathcal{P}}
\newcommand*{\w}{{\mathbf w}}
\def\sR{\mathcal{R}}
\def\sF{\mathcal{F}}
\begin{document}

\title{
On Tackling High-Dimensional Nonconvex Stochastic Optimization via Stochastic First-Order Methods with Non-smooth Proximal Terms and Variance Reduction 
}


\author{Yue Xie     \and Jiawen Bi    \and
        Hongcheng Liu 
}


\institute{Yue Xie \at
              Department of Mathematics and Musketeers Foundation Institute of Data Science, The University of Hong Kong, Pokfulam, Hong Kong.\\
              \email{yxie21@hku.hk}           
           \and
           Jiawen Bi \at
              Department of Industrial and Systems Engineering, The University of Minnesota Twin Cities, Minneapolis, MN.\\
              \email{bjw010529@gmail.com}
           \and
           Hongcheng Liu \at
              Department of Industrial and Systems Engineering, University of Florida, Gainesville, FL. \\
              \email{hliu@ise.ufl.edu} 
}

\date{Received: date / Accepted: date}

\maketitle

\begin{abstract}
When the nonconvex problem is complicated by stochasticity, the sample complexity of stochastic first-order methods may depend linearly on the problem dimension, which is undesirable for large-scale problems. To alleviate this linear dependence, we adopt non-Euclidean settings and propose two choices of non-smooth proximal terms when taking the stochastic gradient steps. This approach leads to stronger convergence metric, incremental computational overhead, and potentially dimension-insensitive sample complexity. We also consider further acceleration through variance reduction which achieves near optimal sample complexity and, to our best knowledge, is the first such result in the $\ell_1/\ell_\infty$ setting. Since the use of non-smooth proximal terms is unconventional, the convergence analysis deviates much from algorithms in Euclidean settings or employing Bregman divergence, providing tools for analyzing other non-Euclidean choices of distance functions. Efficient resolution of the subproblems in various scenarios is also discussed and simulated. We illustrate the dimension-insensitive property of the proposed methods via preliminary numerical experiments. 

\end{abstract}

\section{Introduction}\label{section: intro}
In this paper, we consider a stochastic optimization (SO) problem formulated as below:
\begin{align}\label{SP problem}
\begin{aligned}
        \min & \quad f(\xbf) \textcolor{black}{\triangleq} \bE[F(\xbf,\zeta)] \\
        \mbox{s.t.} & \quad \xbf \in X, 
\end{aligned}
\end{align} 
where $f: \bR^d \to \bR$ and $F: \bR^d \times \Xi \to \bR$ are potentially nonconvex, $(\Xi, \mathcal{F}, \mathbb{P})$ is the probability space, and $X\subseteq \R^d$ is the feasible region.
Throughout the paper, we assume  that $X$ is closed and convex, the expectation   $\E\left[F(\xbf,\zeta)\right] = \int_{\Xi} F(\xbf,\zeta)\,\text{d}\mathbb P(\zeta)$ is  well-defined and finite-valued for every $\xbf\in\R^d$, the expected function $f$ is everywhere differentiable,  it is possible to generate independent and identically distributed realizations, $\zeta_1,\,\zeta_2,\textcolor{black}{\ldots},$ of the random vector $\zeta$ and there exists  a deterministic and measurable function $G: X\times \Xi\rightarrow\textcolor{black}{\R^d}$ such that $\E[G(\xbf,\zeta)]=\nabla f(\xbf)$ for all $\xbf\in X$. As we formalize in Section \ref{sec: Prelim}, some of our results further assume the differentiability of $F(\cdot,\zeta)$, whose gradient is denoted by $\nabla F(\cdot,\zeta)$.  We use $\| \cdot \|$ to denote $\| \cdot \|_2$ in this paper unless stated otherwise.  Problems of this type have been extensively studied and applied in diverse contexts (see, e.g., \cite{lan2020first,shapiro2014lectures,birge1997introduction,ghadimi2013stochastic,fang2018spider,wang2019spiderboost}). In the absence of convexity, our focus follows this line of work by aiming to compute a stationary point with small expected residual, thereby approximating the first-order optimality condition as per Definition~\ref{local stationary point definition}.
\textcolor{black}{A residual function is a nonnegative quantity that measures violation of a first-order stationarity condition. In the constrained setting of \eqref{SP problem}, an exact stationary point satisfies $0\in\partial(f+\delta_X)(\x)$, where $\delta_X$ is the indicator function of $X$ ($\delta_X(\x) = 0$ if $\x \in X$ and $\delta_X(\x) = +\infty$ otherwise), and $\partial$ denotes the Fréchet subdifferential \cite{rockafellar2009variational}.}
\begin{definition}\label{local stationary point definition}
For any $\epsilon\geq 0$ and nonnegative residual function $r(\cdot)$, a feasible random solution $\xbf^{*}_{\epsilon}\in X$ is said to be an {\it $\epsilon$-stationary solution/point} if and only if  
\begin{align}\label{def: epstat}
\bE\left[ r(\xbf^{*}_{\epsilon}) \right]\leq \epsilon.
\end{align}
\end{definition}
It is worth mentioning that Def.~\ref{local stationary point definition} in our analysis is substantiated as an intuitive and natural form of:
\begin{align}\label{def: resfunc}
\begin{aligned}
r(\x) & \triangleq {\rm dist}_{\| \cdot \|_\infty}(0, \partial(f + \delta_X)(\x) ) \triangleq \inf\{ \| \vbf \|_\infty \mid \vbf \in \partial(f + \delta_X)(\x) \}.
\end{aligned}
\end{align}
{\revi The residual specified above is directly tied to problem \eqref{SP problem}: when $r(\x) = 0$, $\x$ is an exact stationary point. } 
\textcolor{black}{Throughout this paper, when we say that one residual is stronger than another, we mean that smallness of the former implies smallness of the latter, up to a dimension-free explicit constant. Thus, convergence in the residual \eqref{def: resfunc} gives at least as sharp an approximation to stationarity as convergence in the alternative residuals compared in Section~\ref{sec: res.f}.}
We would like to point out that, although the residual function is given in $\| \cdot \|_\infty$, the dimension-insensitive result is not trivial. In fact, the residual function is stronger than a residual function discussed in literature (\cite{zhang2018convergence}) and defined in $\| \cdot \|_1$:
\begin{align}
    r_1(\x) \triangleq \frac{1}{\lambda} \| \x - {\prox}_{\lambda f}(\x) \|_1,\label{alternative r}
\end{align}
where $ {\prox}_{\lambda f}(\x) $ is the proximal operator \textcolor{black}{with respect to} Bregman divergence. Detailed discussions are in Section~\ref{sec: res.f}.

One of the mainstream approaches for solving problems  \eqref{SP problem} (so as to attain \eqref{def: epstat})
 is the family of stochastic first-order methods (S-FOMs), \textcolor{black}{also commonly referred to as stochastic approximation methods in the classical literature,} which therefore constitute the focus of this paper. Originating with the seminal work of \cite{robbins1951stochastic,chung1954stochastic,sacks1958asymptotic}, these methods have since been extensively developed and analyzed (see, e.g., \cite{polyak1990new,polyak1992acceleration,nemirovski2009robust,lan2012optimal,ghadimi2013stochastic,rosasco2019convergence,devolder2011stochastic,ghadimi2012optimal,chambolle2018stochastic,chambolle2011first,bach2013non,metel2019simple,xu2019stochastic,nitanda2017stochastic,davis2019stochastic,li2019ssrgd,nguyen2017sarah,wang2019spiderboost,fang2018spider,pham2020proxsarah,horvath2020adaptivity,allen2018make}). These algorithms, including stochastic gradient descent and its numerous variants, are characterized by their reliance on stochastic first-order oracles that yield unbiased estimates of $\nabla f$. 
 
 However, most existing S-FOMs exhibit oracle complexities that scale polynomially with the problem dimension, which becomes prohibitive when $d$ is large. For instance,  in the Euclidean setting, the stochastic gradient descent method yields a complexity of {\revi $\mathcal O(L\Delta_f\sigma^2 \epsilon^{-4})$ stochastic first-order oracle calls, up to constants and lower-order terms, for finding an $\epsilon$-stationary point.  Here $L$ is the Lipschitz smoothness constant for $f$,  $\Delta_f$ is the initial optimality gap, and $\sigma^2$ denotes the variance of the stochastic gradient.  
This complexity holds for either the unconstrained case or the constrained case in the standard smooth nonconvex setting (see \cite{ghadimi2013stochastic} for the unconstrained case, \cite{ghadimi2016mini} for the constrained case and beyond).  Note that $\sigma^2$ may depend on $d$. More specifically, }
\begin{equation*}
\sigma^2 \textcolor{black}{\triangleq} \sup_{\x \in X} \E\left[\Vert G(\x,\zeta) - \nabla f(\x)\Vert^2\right].
\end{equation*}
When $d$ grows, this variance typically scales linearly with $d$, since
\begin{align*}
  \sup_{\x \in X} \E\left[\Vert G(\x,\zeta)-\nabla f(\x) \Vert^2\right] 
= \sup_{\x \in X} \sum_{i=1}^d \E\left[ \big(g_i(\x,\zeta) - \nabla_i f(\x)\big)^2 \right]
= \sO(d),
\end{align*}
where $g_i$ is the $i$th component of $G$ {\revi and $\nabla_i f(\x)$ denotes the $i$th component of $\nabla f(\x)$.}
This rate of dimensional growth directly propagates into the oracle complexity of S-FOMs, leading to a complexity of 
\begin{equation}\label{state-of-the-art-bound}
{\revi \mathcal O\left( \frac{L \Delta_f d}{\epsilon^4}\right),}
\end{equation}
 which is {\revi at least} linear in $d$.  {\revi In fact, $L$ may also scale with $d$. It corresponds to the upper bound on the largest singular value of the Hessian matrix, which could scale when the dimension of the data grows.  For example, consider function $f(\x) \triangleq h(b \cdot a^T \x )$ where $b \in \bR$,  $a \in \bR^d$ and $h$ is the loss function.  The Hessian $\nabla^2f(\x) = h''(b\cdot a^T \x) \cdot b^2 \cdot a a^T $ and its largest singular value $s(\x) = | h''(b\cdot a^T \x) | \cdot b^2 \| a \|^2 $.  Therefore, the estimate of $L$,  which is $\max_\x s(\x) $, scales with $\| a \|^2 = \sO(d)$. }
 
More refined variance-reduction methods \cite{fang2018spider,wang2019spiderboost,zhou2020stochastic} partially improve this dependence and provably achieve $\tilde{\sO}( L \Delta_f \sigma \epsilon^{-3} + \sigma^2 \epsilon^{-2} )$ where $\tilde{\mathcal{O}}$ hides logarithmic factors, which translates to
\begin{align}\label{state-of-the-art-bound epsilon}
{\revi \tilde{\mathcal O}\left(\frac{d}{\epsilon^2}+\frac{L \Delta_f \sqrt{d}}{\epsilon^3}\right)},
\end{align}
if considering the linear dependence of $\sigma^2$ on $d$.  \textcolor{black}{Thus, compared with the basic $\epsilon^{-4}$ dependence in \eqref{state-of-the-art-bound},  \eqref{state-of-the-art-bound epsilon} improves the leading dependence on the target accuracy to $\epsilon^{-3}$ to obtain an $\epsilon$-stationary point with the same residual functions.  This is achieved mainly by incorporating refined variance-reduction techniques, or gradient aggregation,  differentiability of $F(\cdot, \zeta)$ and a stronger assumption on the smoothness condition (mean-square $L$-smoothness or uniform $L$-smoothness across $X \times \Xi$).  The $\epsilon^{-3}$ dependence can be viewed as a benchmark result since a lower bound $const \cdot (L \Delta_f \sigma \epsilon^{-3} + \sigma^2 \epsilon^{-2})$ has been proved in \cite{arjevani2023lower} for the unconstrained case given $\Delta_f$,  $\sigma^2$ and the mean-squared $L$-smoothness condition in the Euclidean setting.}  Nevertheless, the complexity bound remains unfavorable in high dimensions.
 
With the ever-increasing demand for more sophisticated and comprehensive models, stochastic optimization problems of extremely high dimensionality are rapidly emerging. In modern  applications, such as training AI models, it is not uncommon to encounter problems with hundreds of billions of decision variables.
According to  \eqref{state-of-the-art-bound} and \eqref{state-of-the-art-bound epsilon}, SGD and its S-FOM variants may require on the order of hundreds of billions of iterations --- or equivalently, stochastic first-order oracle calls --- to achieve the desired solution accuracy. Without fundamental advances that directly address this curse of dimensionality, solving the SO problem \eqref{SP problem} will become prohibitively difficult in the near future.

Dimension-insensitive S-FOMs have   been made available for convex SO problems  by \cite{nemirovski2009robust,lan2020first}, especially through the non-Euclidean method - stochastic mirror descent (SMD) algorithm which adopts Bregman divergence in the $\ell_1$-norm setting. It has been shown therein that such a method yields a complexity of $$\sO(\max_{\xbf\in X}\E[\Vert G(\xbf,\zeta)\Vert^2_\infty ]/\epsilon^2),$$
in generating a globally $\epsilon$-suboptimal solution. This complexity rate  can often be almost dimension-free, up to a (poly-)logarithmic term, under common conditions. Indeed, it can be verified that $\E[\Vert G(\xbf,\zeta)\Vert^2_\infty ]=\sO(\log d)$ if the gradient estimation errors are sub-Gaussian. In addition, dimension insensitive results are also obtained by  \cite{agarwal2012stochastic} through an successive approximation scheme that solves a series of $\ell_1$-regularized optimization problems using Nesterov’s dual averaging algorithm.  Yet, none of these discussions are readily generalizable to nonconvex problems.

Beyond convex conditions, discussions on dimension-insensitive S-FOMS are relatively scarce.  To our knowledge, perhaps the only prior work in this regard are \cite{zhang2018convergence} and \cite{fatkhullin2024taming}.  The former studies the SMD in solving $\rho$ relatively weakly convex and potentially non-smooth SO problems. For SMD in the $\ell_1$-norm setting, the same work shows a complexity rate of \begin{equation}\label{result historical 2}
{\revi \mathcal O\left(\frac{ \rho \Delta_f \max_{\xbf\in X}\E\left[\Vert G(\xbf,\zeta)\Vert^2_\infty\right]}{\epsilon^4}\right).}
\end{equation}
in achieving an $\epsilon$-stationary point in the sense of \eqref{alternative r} and a so-called Bregman proximal mapping (more discussions are in Section~\ref{sec: Prelim}). In \cite{fatkhullin2024taming}, the authors assumed differentiability of $F(\cdot, \zeta)$,  $\rho$ relatively smoothness of $f$ and \textcolor{black}{derive} the complexity of 
\begin{align}\label{result_fat}
{\revi \mathcal{O}\left( \frac{  \rho \Delta_f \max_{\x \in X} \bE[ \| \nabla F(\x, \zeta) - \nabla f(\x) \|_{\infty}^2 ] }{\epsilon^4} \right)}
\end{align}
in a measure named Bregman Forward-Backward Envelope. \textcolor{black}{Under componentwise sub-Gaussian noise, the quantities $\max_{\xbf\in X}\E[\Vert G(\xbf,\zeta)\Vert^2_\infty]$ in \eqref{result historical 2} and $\max_{\xbf\in X}\E[\Vert \nabla F(\xbf,\zeta)-\nabla f(\xbf)\Vert^2_\infty]$ in \eqref{result_fat} are typically $\mathcal O(\log d)$ rather than linear or polynomial in $d$.  In this sense, these bounds improve the dependence on $d$ in \eqref{state-of-the-art-bound} and \eqref{state-of-the-art-bound epsilon}.} 

In contrast, this paper presents an alternative design  of dimension-insensitive S-FOM, referred to as the DISFOMs which have the following properties:
\begin{itemize}
\item Our algorithm design features non-smooth proximal terms which matches well the non-Euclidean settings including problem assumptions and definition of the residual function, yet seldomly discussed in literature. Under standard assumptions, to obtain an $\epsilon$-stationary point in the sense of \eqref{def: resfunc}, we establish that DISFOM attains the complexity bound
\begin{align}\label{bound1_this_work}
\sO \left(  \frac{ L \Delta_f \max_{\xbf\in X}\E\left[\Vert G(\xbf,\zeta)-\nabla f(\xbf)\Vert_\infty^2\right] \cdot \log d}{\epsilon^4} \right).
\end{align}
Here $L$ is the Lipschitz smoothness constant in the non-Euclidean setting.  {\revi This bound improves the the dependence on dimensionality compared to \eqref{state-of-the-art-bound} and \eqref{state-of-the-art-bound epsilon} under componentwise sub-Gaussian noise assumption.} This bound is comparable to that of SMD in \eqref{result historical 2} and \eqref{result_fat}, which are derived in a more general setting.  However, our result guarantees a sharper approximation to stationarity, since \eqref{def: resfunc} --- which is provably attained by DISFOM --- is a stronger residual measure than the ones used in SMD by \cite{zhang2018convergence} and \cite{fatkhullin2024taming}. Moreover, non-smoothness promotes sparsity - non-smooth proximal terms have the potential to generate sparser solutions, a property relevant in high-dimensional optimization.
\item Furthermore, we prove that a variance-reduced version of DISFOM, referred to as DISFOM\_vr, can achieve a substantial  acceleration through exploiting the common structural assumption that $\nabla F(\cdot,\zeta)$ satisfies the mean-squared $L$-Lipschitz smoothness condition (see more details in  Assumption \ref{ass: Lip} subsequently). The resulting sharpened complexity rate becomes 
\begin{align}\label{bound2_this_work}
 \sO\left( \frac{ L \Delta_f  \max_{\x \in X} \sqrt{ \E[\Vert \nabla F(\xbf,\zeta)-\nabla f(\xbf)\Vert_\infty^2]}\cdot \log^2 d }{\epsilon^3} \right).
\end{align} 
{\revi \eqref{bound2_this_work} further improves over \eqref{result historical 2} and \eqref{result_fat} in terms of dependence on $\epsilon$.  It also improves \eqref{state-of-the-art-bound}, \eqref{state-of-the-art-bound epsilon} in terms of either/both of $d$  and $\epsilon$.  Note that the residual $r(\cdot)$ used to derive \eqref{state-of-the-art-bound} and \eqref{state-of-the-art-bound epsilon} is typically Euclidean norm of the gradient for the unconstrained problems and Euclidean norm of the projected gradient for the constrained problems in literature.  The residual function \eqref{def: resfunc} used in our work is not necessarily stronger than these.  However,  to our best knowledge,  there is no result in literature showing or directly implying that these complexity bounds can be made dimension insensitive ($\log d$ dependence, at least for the variance term) under the componentwise sub-Gaussian noises when the residual function \eqref{def: resfunc} is used and other assumptions including the objective function smoothness remain the same (in this work we have even weaker assumptions). }
\item Last but not least, we discuss the solutions to the subproblems in detail, which is often an overlooked issue for algorithms leveraging non-Euclidean geometry in literature. Despite non-smoothness of the proximal terms, we show that the subproblems can be addressed efficiently and derive the closed-form formula when the problem is unconstrained, or the feasible region are box constraints, or $\ell_1$-ball + box constraints. General linear constraints are also discussed. 
\end{itemize}

\noindent Our numerical experiment verifies efficiency of our subproblem solvers and illustrates the dimension-insensitive property of the proposed algorithms and shows their advantages in  comparison with several other popular algorithms, particularly for solving large-scale SO problems.

We would like to additionally claim a novelty in our algorithm design and analysis. In particular,  the proposed DISFOMs  are based on the use of  non-Euclidean and non-smooth distance metrics to construct per-iteration proximal terms. Consequently, many critical properties of the commonly used Bregmen divergence may no longer be available.  First, the three-point identity --- an underpinning property for the analysis by \cite{zhang2018convergence} and by \cite[Lemmas 6.6 and   3.5, therein]{lan2020first}--- may not hold.  Second, the choice of distance metric may not induced by the inner product, and thus the proof scheme by \cite{ghadimi2013stochastic} would not transfer to our setting, nor would the property of square-integrable martingales --- a seemingly essential element for the analysis by \cite{fang2018spider,wang2019spiderboost} --- easily remain to hold.   Third, there is no known distance generating function  corresponding to  the chosen proximal terms, thus the  Lipschitz property for  the general prox-mapping (as proven by \cite[Lemma 6.5]{lan2020first}) is hardly preserved.  As a result, our proof argument would substantially deviate from many existing theories for S-FOMs. We consider the discussion of such non-standard proximal terms increasingly critical and anticipate the tools derived from our analysis to be growingly useful due to the recent connections drawn (e.g., by \cite{lau2025polargrad}) between the use of non-differentiable proximal terms  and the designs of powerful algorithms, such as ADAM \cite{kingma2014adam} and MUON \cite{jordan2024muon} for AI training.



The rest of the paper is organized as follows: In Section~\ref{sec: Prelim} we discuss our main assumptions and their direct implications. We also discuss the residual function defined in \eqref{def: resfunc} and how it compares with others in literature. In Section~\ref{sec: DISFOM}, we formally introduce the DISFOMs. Theoretical analysis are provided in Section~\ref{sec: analysis}. Section~\ref{sec: subprob} discusses the subproblem solutions. Section~\ref{section:num} includes numerical experiments.
 Section~\ref{sec: conclusion} concludes the paper.
 
\section{Assumptions and discussions of residual function}\label{sec: Prelim}
We consider problem~\eqref{SP problem}  under the following assumptions.


\begin{assumption}\label{assumption first-order oracle}
A stochastic first-order oracle (SFO) is available such that $G(\xbf,\xi)$ can be efficiently evaluated for any given input pair of $(\xbf,\zeta)\in\R^d\times \Xi$.
\end{assumption}
Assumption \ref{assumption first-order oracle} stipulates that the oracle accessible here is the tractable computation of a stochastic gradient. 


\begin{assumption}\label{ass: subG}
    Let $w(\xbf) \triangleq G(\xbf,\zeta) - \nabla f(\xbf)$. 
    There exists some $\sigma_\infty^2 < \infty$ such that $\max_{\x \in X} \bE[ \| w(\x) \|_\infty^2 ] = \sigma_\infty^2$.
\end{assumption}
Assumption~\ref{ass: subG} assumes that besides unbiasedness, the stochastic gradient error has bounded variance in $\ell_\infty$-norm. If we suppose that $w_i(\x)$ conforms sub-Gaussian distribution with sub-Gaussian norm $\sigma$ across all the indices and $X$, then $\sigma_\infty^2 = \mathcal{O}( (\log d) \sigma^2)$. See Lemma~\ref{lm: subG} in Appendix~\ref{app: proof} for more details.
\begin{assumption}\label{asp: flip} There exists some constant $L\geq 0$ such that
\begin{align}\label{ineq: lip}
        \| \nabla f(\x) - \nabla f(\y) \|_\infty \le L \| \x - \y \|_1, \quad \forall \x, \y \in X.
\end{align}
\end{assumption}
In particular, Assumption~\ref{asp: flip} indicates that the following holds:
\begin{align}\label{ineq: Taylor}
    f(\ybf) \le f(\xbf) + \nabla f(\xbf)^T (\ybf - \xbf) + \frac{L}{2} \| \ybf - \xbf \|_1^2.
\end{align}
Assumption~\ref{asp: flip} is necessary in the first series of analysis of our proposed algorithm.  To further improve the complexity guarantees, we need differentiability of $\nabla F(\cdot, \zeta)$ and a stronger assumption regarding smoothness of the objective function as follows.
\begin{assumption}[Mean-squared Lipschitz smoothness condition]\label{ass: Lip} \phantom{.}\\
Function $F(\cdot,\zeta)$ is everywhere differentiable for all $\zeta\in \Xi$. Its gradient $\nabla F(\cdot,\zeta)$ satisfies that, for some constant $L\geq 0$,  
    \begin{align}
    \begin{aligned}
        \bE[\| \nabla F(\xbf,\zeta) - \nabla F(\ybf,\zeta) \|_\infty^2] \le L^2 \| \xbf - \ybf \|_1^2, \quad \forall \xbf, \ybf \in X.
        \end{aligned}
    \end{align}
\end{assumption}
Assumption~\ref{ass: Lip} implies Assumption~\ref{asp: flip} by Jensen's inequality. It is vital to the analysis when incorporating variance reduction in the algorithm. Last, we need a mild assumption regarding the lower boundedness of the objective function.
\begin{assumption}\label{ass: setw3hj} There exists some   $\lf\in\R$ such that
\begin{align}
    f(\xbf) \ge \lf, \quad \forall \xbf \in X.
\end{align}
\end{assumption}

All these assumptions above together with the basic assumptions in Section~\ref{section: intro} are standard or equivalent to the common conditions (or their counterparts in non-Euclidean norms) in the literature of stochastic approximation or S-FOMs. For example, Assumptions \ref{assumption first-order oracle} are imposed by \cite{nemirovski2009robust}.  Assumption~\ref{asp: flip} (Assumption \ref{ass: Lip}) in Euclidean norm is necessary for the discussions by \cite{ghadimi2013stochastic,ghadimi2016mini} (\cite{fang2018spider,wang2019spiderboost,nguyen2017sarah,pham2020proxsarah}).
{\revi In fact, Assumption~\ref{asp: flip} and \ref{ass: Lip} are weaker than their counterparts in the Euclidean setting.  More discussions can be found in Remark~\ref{rmk: Lip.const.} and Appendix~\ref{app: Lip}. } 

\subsection{Discussion of residual function}\label{sec: res.f}
A key ingredient of the worst-case complexity guarantees is the criterion used to measure optimality.  In this subsection we do not enforce blanket assumptions on the function $f$ other than differentiability over an open region including $X$.  We demonstrate that the natural definition \eqref{def: resfunc} of residual function is stronger than many existing measures in literature including:
\begin{align*}
    & r_1(\x) \triangleq \frac{1}{\lambda} \| \x - {\rm prox}_{\lambda f}(\x) \|_1, \\
    & r_2(\x) \triangleq \frac{1}{\lambda} \sqrt{ D_h(\x,{\prox}_{\lambda f}(\x)) + D_h({\prox}_{\lambda f}(\x),\x)},\\
    & r_3(\x) \triangleq \sqrt{-\min_{\y \in X} Q(\x,\y)},
\end{align*}
where 
\begin{align}\label{def: prox}
{\rm prox}_{\lambda f}(\x) \triangleq \mbox{arg}\min\limits_{\z \in X}{ f(\z) + \frac{1}{\lambda} D_h(\z,\x) },
\end{align}
and $D_h(\cdot,\cdot)$ is the Bregman divergence induced by \textcolor{black}{a} 1-\textcolor{black}{strongly} convex continuously differentiable \textcolor{black}{distance-generating} function $h(\cdot)$ \textcolor{black}{with respect to} $\| \cdot \|_1$, i.e.,
\[
D_h(\ubf,\vbf) \triangleq h(\ubf) - h(\vbf) - \langle \nabla h(\vbf), \ubf - \vbf \rangle \ge \frac{1}{2}\| \ubf - \vbf \|_1^2,
\]
and
\[
Q(\x,\y) \triangleq \langle \nabla f(\x), \y - \x \rangle + \frac{1}{\lambda} D_h(\y, \x).
\]
{\revi The residuals $r_1(\x)$ and $r_2(\x)$ are the residual measures used in the mirror descent analysis of \cite{zhang2018convergence}; $r_1(\x)$ is also the $\ell_1$-norm analogue of the residual used in \cite{davis2018stochastic}.  $r_2(\x)$ is named Bregman Proximal Mapping. The residual $r_3(\x)$ is the Bregman Forward Backward Envelope residual used in \cite{fatkhullin2024taming}.  When we say that one residual is stronger than another, we mean that smallness of the former implies smallness of the latter, up to a dimension-free explicit constant. } In fact, $r_2(\x)$ is stronger than $r_1(\x)$ since 
\begin{align*}
    & D_h(\x,{\prox}(\x)) + D_h({\prox}(\x),\x) = \langle \nabla h(\x) - \nabla h({\prox}(\x)), \x - {\prox}(\x) \rangle \\
    & \ge \| \x - {\prox}(\x) \|_1^2.
\end{align*}

First, we point out that $r(\x)$ is stronger than $r_2(\x)$. In fact, if $f$ is in addition $\rho$-relatively weakly convex (i.e., $f(\x) + \rho h(\x)$ is convex for such $h(\cdot)$ and some $\rho > 0$, a weaker condition than Assumption~\ref{asp: flip} when $\rho = L$) and $\lambda \in (0,\rho^{-1})$ so that \eqref{def: prox} is well-defined,  we have the following statement.
\begin{lemma}\label{lm: comp.mea}
Suppose that $f$ is $\rho$-relatively weakly convex \textcolor{black}{with respect to} $h$ \textcolor{black}{and that $\lambda \in (0,\rho^{-1})$}. Let $\x \in X$, then
    \begin{align}\label{ineq: rf0}
    \begin{aligned}
 r(\x) \ge (1-\lambda \rho)r_2(\x).
 \end{aligned}
\end{align}
\end{lemma}
\begin{proof}
    Denote $\hx \triangleq {\prox}_{\lambda f}(\x)$. Then by definition \eqref{def: prox}, 
\begin{align}
\notag
    & 0 \in \nabla f(\hx) + \frac{1}{\lambda}( \nabla h(\hx) - \nabla h(\x) ) + \sN_X(\hx), \\
    \notag
    & \Longleftrightarrow \left\langle \left( \nabla f(\hx) + \frac{1}{\lambda} ( \nabla h(\hx) - \nabla h(\x) ) \right) , \y - \hx \right\rangle \ge 0, \quad \forall \y \in X, \\
    \notag
    & \implies \langle \nabla f(\x), \x - \hx \rangle \ge \left\langle \nabla f (\hx) + \frac{1}{\lambda} \nabla h(\hx) - \nabla f(\x) - \frac{1}{\lambda} \nabla h(\x), \hx - \x \right\rangle \\
    \label{ineq: rf1}
    & \qquad \ge \left( \frac{1}{\lambda} - \rho \right) \langle \nabla h(\hx) - \nabla h(\x), \hx - \x \rangle.
\end{align}
Choose arbitary $\vbf \in \partial(f + \delta_X)(\x)$.  Then the differentiability of $f$ and convexity of $\delta_X(\x)$ imply that $\vbf \in \nabla f (\x) + \partial \delta_X(\x) = \nabla f (\x) + \sN_X(\x)$ and we have
\begin{align}\label{ineq: rf2}
    \left\langle \vbf - \nabla f(\x), \y - \x \right\rangle \le 0, \quad \forall \y \in X.
\end{align}
By \eqref{ineq: rf1} and \eqref{ineq: rf2}, we have
\begin{align}\label{ineq: rf3}
    \left( 1/\lambda - \rho \right) \langle \nabla h(\hx) - \nabla h(\x), \hx - \x \rangle \le \langle \vbf , \x - \hx \rangle.
\end{align}
\textcolor{black}{Without loss of generality}, suppose that $\hx \neq \x$. Then by \eqref{ineq: rf3},
\begin{align}\label{ineq: rf4}
    \left( 1/\lambda - \rho \right) \frac{\langle \nabla h(\hx) - \nabla h(\x), \hx - \x \rangle} {\| \hx - \x \|_1} \le \left\langle \vbf , \frac{\x - \hx}{\| \x - \hx \|_1} \right\rangle \le \| \vbf \|_\infty.
\end{align}
However, 
\begin{align}
\notag
    & \langle \nabla h(\hx) - \nabla h(\x), \hx - \x \rangle \ge \| \hx - \x \|_1^2 \\
    \label{ineq: rf5}
   \implies &  \frac{\langle \nabla h(\hx) - \nabla h(\x), \hx - \x \rangle} {\| \hx - \x \|_1} \ge \sqrt{\langle \nabla h(\hx) - \nabla h(\x), \hx - \x \rangle}.
\end{align}
By combining \eqref{ineq: rf4}\eqref{ineq: rf5} and the fact that $\vbf$ is arbitrary, we have that \eqref{ineq: rf0} holds. \qed
\end{proof}

Second, we illustrate that $r(\x)$ is stronger than $r_3(\x)$ via the following Lemma.
\begin{lemma}
Suppose that $\x \in X$. Then
\begin{align*}
    r(\x) \ge \sqrt{ 2 / \lambda }  \cdot r_3(\x).
\end{align*}
\end{lemma}
\begin{proof}
    Let $\hat \y \triangleq \mbox{argmin}_{\y \in X} Q(\x,\y)$. Then 
    \begin{align*}
        r_3(\x) = \sqrt{ \langle \nabla f(\x), \x - \hat \y \rangle - (1/\lambda)D_h (\hat \y, \x) }.
    \end{align*}
    Choose arbitary $\vbf \in \partial(f + \delta_X)(\x)$, then we have
\begin{align}\label{ineq: rf2-2}
    \left\langle \vbf - \nabla f(\x), \y - \x \right\rangle \le 0, \quad \forall \y \in X.
\end{align}
Therefore, by pluing $\y = \hat \y$ in \eqref{ineq: rf2-2}, we have
\begin{align*}
    \langle \nabla f(\x), \x - \hat \y \rangle & \le \langle \vbf, \x - \hat \y \rangle \le \| \vbf \|_\infty \| \x - \hat \y \|_1 \le \frac{\lambda}{2} \| \vbf \|_\infty^2 + \frac{1}{2 \lambda} \| \x - \hat \y \|_1^2 \\
    & \le \frac{\lambda}{2} \| \vbf \|_\infty^2 + \frac{1}{\lambda} D_h( \hat \y, \x ) \implies r_3^2(\x) \le (\lambda/2)\| \vbf \|_{\infty}^2.
\end{align*}
Note that $\vbf$ is arbitrary, so we conclude the result.\qed
\end{proof}
{\revi 
In fact,  $r(\x)$ is less sensitive to, or rather not affected by, the Lipschitz smoothness constant of the distance generating function $h$ compared to other residuals.  To leverage the dimension-insensitive effect of non-Euclidean methods,  one needs to let the distance generating function $h$ be strongly convex with respect to the  $\ell_1$-norm in those Bregman/mirror-descent approaches. The typical choices are $h(\x) \triangleq \sum_{i=1}^d \x_i \log (\x_i)$ and $h(\x) \triangleq \frac{1}{2} \| \x \|_p^2$ ($p > 1$, $p \approx 1$).  Note that these functions are potentially ill-conditioned near $0$ (or when some coordinates are near $0$).  Their Lipschitz smoothness constants blow up near $0$,  which significantly reduces the strength of the convergence metric in those works.  Let us consider comparison with the Bregman proximal mapping $r_2 (\x)$ 
and the Bregman Forward Backward Envelope $r_3 (\x)$.
Let
\begin{align*}
T_{\lambda, f}(\x) \triangleq  {\rm arg}\min\limits_{ \y \in X} Q(\x,\y) = {\rm arg}\min\limits_{ \y \in X} \left\langle \nabla f(\x), \y - \x \right\rangle + \frac{1}{\lambda} D_h( \y,\x).
\end{align*}
Then $r_3(\x) = \sqrt{-Q(\x,T_{\lambda,f}(\x))}$.  Consider the following example:
\begin{align*}
& \nabla f(\x) \equiv (1,0,...,0)^T \triangleq {\bf e}_1, \; X = \{ \x \ge 0 \mid \| \x \|_1 \le 1 \}.  
\end{align*}
Suppose that $\tilde{\x} = (\epsilon,1/d,...,1/d) $ and $\epsilon \in (0,1/d)$ is small enough.  Then $r(\tilde \x) = \| \nabla f(\tilde \x) \|_\infty \equiv 1$. Denote $\zbf \triangleq {\rm prox}_{\lambda f}(\tilde \x) $ and 
$\y \triangleq T_{\lambda,f}(\tilde \x) $.  Then by optimality conditions and the fact that $\nabla f$ is a constant,  we have that $\y = \zbf$ (we will use them interchangeably) and 
\begin{align}\label{hoptcon}
( \nabla h(\y) - (\nabla h(\tilde \x) - \lambda \cdot {\bf e}_1 ) )^T (\x - \y) \ge 0, \quad \forall \x \in X.
\end{align}
Therefore, 
\begin{align*}
\lambda^2 r_2^2(\tilde \x) 
& = (\nabla h(\tilde \x)- \nabla h(\y))^T(\tilde \x - \y)\\
& = (\nabla h(\tilde \x)- \nabla h(\y))^T ( \nabla^2 h(\y + \theta(\tilde \x- \y)) )^{-1}  (\nabla h(\tilde \x)- \nabla h(\y)), \\
r_3^2(\tilde \x) & = \left\langle \nabla f(\tilde \x), \tilde \x - \y \right\rangle - (1/\lambda) D_h( \y,\tilde \x) \le  \left\langle \nabla f(\tilde \x), \tilde \x - \y \right\rangle \\
& = ( \nabla f(\tilde \x) )^T ( \nabla^2 h(\y + \theta(\tilde \x- \y)) )^{-1}  (\nabla h(\tilde \x)- \nabla h(\y)),
\end{align*}
for some $\theta \in [0,1]$.  If $ \nabla^2 h(\y + \theta(\tilde \x- \y))$ is unbounded, i.e., the Lipschitz smoothness constant of $h$ is unbounded, then the above quantities are potentially small.  More specifically,  if we let $h(\x) \triangleq \sum_{i=1}^d \x_i \log \x_i$, then we can deduce that $\nabla h(\y) = \nabla h(\tilde \x) - \lambda \cdot {\bf e}_1$,  $\y_1 = \epsilon \exp(-\lambda) < \tilde{\x}_1$,  $\y_i = \tilde{\x}_i$ for $i = 2,...,d$,
\begin{align*}
\lambda^2 r_2^2(\tilde \x) & = \lambda^2 ( \nabla^2 h(\y + \theta(\tilde \x- \y)) )^{-1}_{1,1} = \lambda^2 (\y_1 + \theta( \tilde{\x}_1) - \y_1)  \le \lambda^2 \tilde \x_1 \\
& = \lambda^2 \epsilon \to 0,  \\
r_3^2(\tilde \x) & = \left\langle \nabla f(\tilde \x), \tilde \x - \y \right\rangle - (1/\lambda) D_h( \y,\tilde \x) \le  \left\langle \nabla f(\tilde \x), \tilde \x - \y \right\rangle = \tilde \x_1 - \y_1 \\
& = (1- \exp(-\lambda)) \epsilon \to 0,
\end{align*}
as $\epsilon \to 0$. If we let $h(\x) \triangleq \frac{1}{2} \| \x \|_p^2$, $p \in (1,2)$, then $[ \nabla h(\x) ]_i = \| \x \|_p^{2-p} | \x_i |^{p-1} \sgn (\x_i)$ and when $\epsilon$ is small, 
we have $[ \nabla h(\tilde \x) ]_1 < \lambda$. Therefore, we deduce that $\y_1 = 0$,  and $[\nabla h(\y)]_i = [ \nabla h(\tilde \x) ]_i$ for $i = 2,...,d$ given small enough $\epsilon$. In such a case,
\begin{align*}
\lambda^2 r_2^2(\tilde \x) & =  [ \nabla h(\tilde \x)]_1 \tilde \x_1 = \| \x \|_p^{2-p} \x_1^p \le \| \x \|_1^{2-p} \x_1^p  \le 1 \cdot \epsilon^p  =  \epsilon^p \to 0, \\
r_3^2(\tilde \x) & = \left\langle \nabla f(\tilde \x), \tilde \x - \y \right\rangle - (1/\lambda) D_h( \y,\tilde \x) \le  \left\langle \nabla f(\tilde \x), \tilde \x - \y \right\rangle = \tilde \x_1 - \y_1 \\
& = \epsilon \to 0,
\end{align*}
as $\epsilon \to 0$.  In either case above, we can see that the ratios $r(\x) / r_2(\x)$ and $r(\x) / r_3(\x)$ (therefore $r(\x) / r_1(\x)$) can be arbitrarily large,  indicating the strength of $r(\x)$ as a convergence metric.
} 

To conclude this section, we give some notations to be used in this work.
\paragraph{Notation.} Given a sequence $\{ \vbf^k \}$, $\Delta \vbf^{k+1} \triangleq \vbf^{k+1} - \vbf^k$. Given a vector $\vbf \in \bR^d$, $\vbf_i$ means its $i$th component unless specified otherwise in the context. For a positive integer $K$, $[K]$ denotes the set $\{1,2,\hdots,K\}$. $\| \cdot \|$ denotes the $2$-norm $\| \cdot \|_2$. Denote $g^k \triangleq \nabla f(\xbf^k)$ for any $k$. For a scalar $\lambda$, define $\sgn(\alpha) \textcolor{black}{\triangleq} \begin{cases} 1 ,& \alpha > 0 \\-1, & \alpha <0  \\ 0, & \alpha = 0 \end{cases}$. Suppose $\operatorname{clip}_{[l_i, u_i]}(\alpha) \triangleq \begin{cases}
    l_i, & \alpha \le l_i \\ \alpha, & \alpha \in (l_i,u_i) \\ u_i, & \alpha \ge u_i
\end{cases}. $ For vectors $\vbf,l,u$, we define $\operatorname{clip}_{[l,u]}(\vbf)$ analogously in a componentwise fashion. ${\rm ri}(S)$ denotes relative interior point of a set $S$. For a proper function $h$, its domain is denoted as ${\rm dom}\; h \triangleq \{ \x \mid h(\x) < + \infty \}$.  We use $\partial$ to denote the Fréchet subdifferential,  while in convex cases it coincides with the classic subdifferential defined for convex functions.

\section{Dimension-insensitive stochastic first-order methods}\label{sec: DISFOM}

\textcolor{black}{This section first presents a generic stochastic first-order template, then specifies the choices for non-smooth proximal terms and stochastic-gradient estimators that lead to dimension-insensitive complexity bounds. The phrase ``dimension-insensitive'' refers to complexity bounds whose dependence on $d$ coming from the variance term can be logarithmic, or polylogarithmic, under the componentwise sub-Gaussian noises.  The bounds therefore may still contain $\log d$ terms, but avoid the linear dependence that arises from Euclidean variance measures. }

We start from a general framework of a stochastic first-order method as follows.

\begin{algorithm}[H] 
\caption{Stochastic first-order method}\label{alg: DI-SGD}
\begin{description}
\item[{\bf Step 1.}] Initialize $\xbf^1 \in X$ and hyper-parameters $K$, $\eta_k > 0$ for all $k = 1,...,K$.
\item[{\bf Step 2.}] Invoke the following operations for $k=1,...,K$:
\begin{description}
 \item[{\bf Step 2.1}] Generate the gradient estimate $G^k$. 
 \item[{\bf Step 2.2}] Solve the proximal projection problem:
 \begin{align}\label{proj}
     \xbf^{k+1} {\revi \triangleq} P_X^k (\xbf^k - \eta_k G^k)
 \end{align}
\end{description}
\item[{\bf Step 3.}] Output $\xbf^{Y+1}$ for a random $Y$, which has a discrete distribution on $[K]$ with a probability mass function $\mathbb P[Y=k]= p_k$.
\end{description}
\end{algorithm}

Now we discuss the details of the proximal projection operator $P_X^k(\cdot)$ and gradient estimate $G^k$ that can achieve the dimension-insensitive property. In particular, the proximal projection operator $P_X^k(\cdot)$ is defined as
\begin{align}\label{def: proj}
    P_X^k(\vbf) \textcolor{black}{\triangleq} \mbox{arg}\min_{\xbf \in X} \left\{ \frac{1}{2} \| \xbf - \vbf \|^2 + \phi(\xbf - \x
^k) \right\},
\end{align}
where $\phi: \bR^d \to \bR \cup \{ + \infty \}$ is a proper closed convex function (possibly non-smooth) and we suppose that calculation of \eqref{def: proj} is exact. We will discuss the efficient solution to \eqref{def: proj} in Section~\ref{sec: subprob}. Next we specify our choices of $\phi$ and $G^k$.

\subsection{Choices of $\phi$}

We provide two possible choices of $\phi$ in \eqref{def: proj}. \\
\noindent{Case 1}: Let 
\begin{align}\label{setting1}
\phi(\x) \textcolor{black}{\triangleq} \frac{\hrho}{2} \| \x \|_1^2, \quad \hrho > 0.
\end{align}
For such a $\phi(\x)$, the following lemma holds by application of Danskin's Theorem.
\begin{lemma}\label{lm: l1subdiff}
Let $\phi(\x)$ be in \eqref{setting1}. 
Then
\begin{align}\label{l1subdiff}
    \partial \phi(\x) & = \hrho \| \x \|_1 \partial \| \x \|_1 = \left\{ \hrho \| \x \|_1 \vbf: {\revi \mbox{for any $i$}}, \begin{array}{ll}
         \vbf_i = 1 & \mbox{if } \x_i > 0 \\
         \vbf_i = -1 & \mbox{if } \x_i < 0\\
         \vbf_i \in [-1,1] & \mbox{if } \x_i = 0
    \end{array} \right\}.
\end{align}
\end{lemma}
Therefore, it is easy to see that $\phi(\x)$ is non-smooth. It is not strongly convex either. Indeed, we could solve the following problem ($\vbf \neq 0$):
\begin{align*}
    \min_{\x} \quad \vbf^T \x + \frac{1}{2} \| \x \|_1^2.
\end{align*}
Based on the optimality conditions: $-\vbf \in \| \x^* \|_1 \partial \| \x^* \|_1$, we have $\| \x^* \|_1 = \| \vbf \|_\infty$; $\x^*_i = 0$ if $|\vbf_i| < \| \vbf \|_\infty$; $\x_i^* = 0$ or $ \sgn(\x^*_i) = - \sgn(\vbf_i) $ if $|\vbf_i| = \| \vbf \|_\infty$. Therefore, if $i \neq j$ and $|\vbf_i| = |\vbf_j| = |\vbf|_\infty$, there will be multiple solutions, negating strong convexity of $\| \cdot \|_1^2$. Therefore, it is impossible to represent $\| \y - \x \|_1^2$ in the form of $ h(\y) - h(\x) - \langle \xi_\x, \y - \x \rangle$, $\xi_\x \in \partial h(\x)$ for some convex $h(\cdot)$. On the other hand, this choice of $\phi$ not only enables convergence (c.f. Section~\ref{sec: analysis}) but also leads to closed-form solutions to the subproblem \eqref{def: proj} for many scenarios including the unconstrained case, $X$ is box-constraints and/or $\ell_1$-ball  (c.f. Section~\ref{sec: subprob}).  An alternative to the choice of $\phi$ in case 1 is an indicator function defined as follows.\\
Case 2: Let 
\begin{align}\label{case2phi}
    \phi(\z) \textcolor{black}{\triangleq} \delta_{\bar X} (\z) \textcolor{black}{\triangleq} \begin{cases}
        0 & \mbox{if } \z \in \bar X \\
        +\infty & \mbox{if } \z \notin \bar X,
    \end{cases}
\end{align}
where $ \bar X \triangleq \{ \z \mid \| \z \|_1 \le \psi \}$ and $\psi > 0$. The non-smooth and non-strongly convex nature of this choice of $\phi$ is more obvious. However, it also facilitates convergence and computes efficiently. 

\subsection{Choices of $G^k$}
Next we discuss choices of the gradient estimate $G^k$ and its variance in the non-Euclidean norm.\\
{\bf Minibatch.} First we consider using a finite batch of samples to estimate the gradient $g^k \triangleq \nabla f(\x^k)$. Although the batchsize is not necessarily small, we borrow the name \textit{minibatch} considering that accurate estimation of $\nabla f(\x^k)$ will take probably huge or even infinite number of samples. More specifically, let
\begin{align}\label{def: Gk-SGD}
    G^k \triangleq \frac{1}{m_k} \sum_{i = 1}^{m_k} G(\x^k,\zeta^k_i),
\end{align}
where $\zeta^k_i$ are i.i.d realizations of $\zeta$. Then the following results hold. 
\begin{lemma}\label{var-minibatch}
    Assumption~\ref{ass: subG} holds and consider the minibatch sampling \eqref{def: Gk-SGD}. Let $\bar w^k \triangleq \frac{1}{m_k} \sum_{i=1}^{m_k} (G(\x^k, \zeta^k_i) - \nabla f(\x^k))$, then for any $k$,
    \begin{itemize}
        \item[1.] $ \bE[\bar w^k \mid \x^k] = 0 $,
        \item[2.] $ \bE[ \| \bar w^k \|_\infty^2 \mid \x^k] \le {\revi e^2} (\log d) \sigma_\infty^2/m_k $.
    \end{itemize}
\end{lemma}
\begin{proof}
\textcolor{black}{The proof uses the following auxiliary function construction. } Let $\| \cdot \|_p$ and $\| \cdot \|_q$ denote the $p$- and $q$- norm of vectors respectively, with $1/p + 1/q = 1$, $p \in [1,+\infty]$, $q \in [1,+\infty]$. $W$ denotes a smooth and 1-strongly convex function \textcolor{black}{with respect to} $\| \cdot \|_p$. Also suppose that $W(x) \le (\Omega/2) \| x \|_p^2 $ and $W(x) \ge 0$. An example of $W$ when $p = 1$ is 
\begin{align}\label{def: W}
W(x) \triangleq (1/2) e (\log d) d^{(\bar p - 1)(2 - \bar p)/\bar p} \| x \|_{\bar p}^2 ,
\end{align}
where $\bar p \textcolor{black}{\triangleq} 1 + 1/(\log d)$. In this case, $\Omega \textcolor{black}{\triangleq} e^2 \log d$ since $\| x \|_{\bar p} \le \| x \|_1$. Such $W(x)$ is $1$-strongly convex \textcolor{black}{with respect to} $\| \cdot \|_1$ since $\| x \|_1 \le \| x \|_{\bar p} d^{\frac{\bar p - 1}{\bar p}} $ (H{\"o}lder's inequality) and $(1/2)\| \cdot \|_p^2$ is $(p-1)$-strongly convex \textcolor{black}{with respect to} $\| \cdot \|_p$ for any $p \in (1,2]$ according to \cite[Lemma 8.1]{ben2001ordered},\cite{nemirovskij1983problem}. When $p \in (1,2]$, we could let $W(x) \triangleq \frac{1}{2(p-1)} \| \cdot \|_p^2 $. In this case, $\Omega \textcolor{black}{\triangleq} 1/(p-1)$.

Now define the dual function as the Legendre transformation $W^*(z) \triangleq \{ \max_x \langle z, x \rangle - W(x) \}$. Then $W^*(0) = 0$. We can also show that
\begin{align*}
    W^*(z) \ge \max_x \{ \langle z,x \rangle - (\Omega/2) \| x \|_p^2 \} = (1/(2\Omega)) \| z \|_q^2,
\end{align*}
where the last equality is because
\begin{align*}
    & \max_x \left\{ \langle z,x \rangle - (\Omega/2) \| x \|_p^2 \right\} \\
    & \le \max_x \left\{ \| z \|_q \| x \|_p - (\Omega/2) \| x \|_p^2 \right\} \\
    & = \max_x \left\{ - (\Omega/2) ( \| x \|_p^2- (2/\Omega) \| z \|_q \| x \|_p + (1/\Omega^2) \| z \|_q^2  ) + (1/(2\Omega)) \| z \|_q^2 \right\} \\
    & = (1/(2\Omega)) \| z \|_q^2,
\end{align*}
and the ``$=$" holds when $\langle z, x \rangle = \| x \|_p \| z \|_q $ and $\| x \|_p = (1/\Omega) \| z \|_q$. Therefore, $W^*(z) \ge 0$ and $\nabla W^*(0) = 0$. Note that $W^*$ is continuously differentiable by \cite[Proposition 12.60]{rockafellar2009variational} since both $W$ and $W^*$ are proper, lsc, convex function. Moreover, $W^*$ is $1$-Lipschitz smooth \textcolor{black}{with respect to} $\| \cdot \|_q$. Indeed, for any $\x, \y, \z_x, \z_y$ such that $\z_x = \nabla W^*(\x)$ and $\z_y = \nabla W^*(\y)$, we have $\x = \nabla W(\z_x)$ and $\y = \nabla W(\z_y)$. Therefore,
\begin{align*}
& \langle \nabla W^*(\x) - \nabla W^*(\y), \x - \y \rangle = \langle \z_x - \z_y, \nabla W(\z_x) - \nabla W(\z_y) \rangle \\
& \ge \| \z_x - \z_y \|_p^2 = \| \nabla W^*(\x) - \nabla W^*(\y) \|_p^2 \\
\implies & \| \nabla W^*(\x) - \nabla W^*(\y) \|_p \le \| \x - \y \|_q.
\end{align*}
Item 1 is immediate from Assumption~\ref{ass: subG}. Now we prove 2. Denote 
    $$
    w^{k,i} \triangleq G(\x^k, \zeta^k_i) - \nabla f(\x^k).
    $$
   {\revi \textcolor{black}{Without loss of generality,} suppose that $m = m_k$. } Then
    \begin{align*}
        {\bar w}^k = \frac{1}{m}\sum_{i=1}^{m} w^{k,i}.
    \end{align*}
    Therefore, by 1-Lipschitz smoothness of $W^*$, we have that
    \begin{align*}
    & W^*( \bar w^{k}) \\
    & \le W^*\left( \frac{1}{m}\sum_{i=1}^{m-1} w^{k,i} \right) + \left\langle \nabla W^*\left( \frac{1}{m}\sum_{i=1}^{m-1} w^{k,i} \right), \frac{1}{m} w^{k,m} \right\rangle + \frac{1}{2m^2} \| w^{k,m} \|_q^2.
    \end{align*}
    {\revi Define
    \begin{align*}
    S_j^k \triangleq \frac{1}{m}\sum_{i=1}^j w^{k,i},  \qquad S_0^k \triangleq 0,
    \end{align*}
    For each $j=1,\ldots,m$, the same smoothness argument gives
    \begin{align*}
    W^*\left( S_j^k \right) & \le W^*\left( S_{j-1}^k  \right)+
    \left\langle \nabla W^*\left( S_{j-1}^k  \right),\frac{1}{m}w^{k,j}\right\rangle  +  \frac{1}{2m^2}\|w^{k,j}\|_q^2.
    \end{align*}
    Conditioning on $\x^k,\zeta^k_1,\ldots,\zeta^k_{j-1}$, the inner-product term has mean zero. Summing these inequalities for $j=2,\ldots,m$ and then taking conditional expectation given $\x^k$ yields}
    \begin{align}\label{ineq1: mininf}
        \bE[ W^*( \bar w^k ) \mid \x^k ]
        \le \frac{1}{2m^2} \sum_{i=1}^{m} \bE[\| w^{k,i} \|_q^2 \mid \x^k].
\end{align}
Let $p\textcolor{black}{\triangleq}1$, $q \textcolor{black}{\triangleq} +\infty$, and choose $W$ as in \eqref{def: W}, we have
\begin{align*}
    \frac{1}{2\Omega} \bE[ \| \bar w^k \|_\infty^2 \mid \x^k ] \le \bE[ W^*( \bar w^k ) \mid \x^k ] \overset{\eqref{ineq1: mininf}}{\le} \frac{1}{2m^2} \sum_{i=1}^{m} \bE[\| w^{k,i} \|_\infty^2 \mid \x^k] \le \frac{\sigmainf^2}{2m},
\end{align*}
{\revi where the last inequality is from Assumption~\ref{ass: subG}. } Plug in $\Omega = e^2 \log d$ and we have the result. \qed
\end{proof}
{\revi
\begin{remark}
Infinity norm square of summands is unlike the 2-norm square where the inner products are canceled when taking (conditional) expectation, making it more difficult to show reduction of variance by sample averaging. Hence an auxiliary function $W$ is needed.  Construction of the auxiliary function is inspired by \cite[Lemma 1]{ilandarideva2024accelerated}. 
\end{remark}
}

\noindent{\bf Variance reduction.} We make use of the variance reduction technique introduced in the SPIDER algorithm \cite{fang2018spider} for further acceleration. The technique requires using relatively large number of samples intermittently when estimating $g^k$, and reuse this relatively more accurate estimation in other iterations to improve sampling efficiency. For the rest of this subsection we assume differentiability of $F(\cdot, \zeta)$.

Given interval length $q_0$, define $n_k \triangleq \lfloor \frac{k-1}{q_0} \rfloor q_0 + 1$. Then $G^k$ is computed as follows:
\begin{align}\label{def: Gk-SPIDER}
& G^k \triangleq \\
\notag
& \begin{cases}
    \frac{1}{m_k} \sum_{i=1}^{m_k} \nabla F(\x^k, \zeta^k_i),
    & \mbox{if } \textcolor{black}{k \equiv 1 \pmod {q_0}},\\
    \frac{1}{m_k} \left( \sum_{i=1}^{m_k} \nabla F(\xbf^k,\zeta^k_i) - \sum_{i=1}^{m_k} \nabla F(\xbf^{k-1},\zeta^k_i) \right) + G^{k-1},  & \mbox{otherwise,}
\end{cases}
\end{align}
where $\zeta^k_i$ are i.i.d realizations of $\zeta$. Similar to Lemma~\ref{var-minibatch}, for any $n_k$,
\begin{align}\label{redvar-SPIDER1}
\bE [ \| g^{n_k} - G^{n_k} \|_\infty^2 ] \le \frac{{\revi e^2}(\log d)\sigma_\infty^2}{m_{n_k}}.    
\end{align}
Suppose that 
\begin{align}\label{mk_spider}
\begin{cases}
    m_k \equiv m & \mbox{if } k \neq n_k \\
m_k \equiv m_1  & \mbox{otherwise}.
\end{cases}
\end{align} 
Likewise, we face the challenge of non-Euclidean norms incompatible with the natural inner-product. First we clarify a notation. Suppose that
$$
B^k(\,\cdot\,)\textcolor{black}{\triangleq}\frac{1}{m_k}  \sum_{i=1}^{m_k} \nabla F(\,\cdot\,,\zeta^k_i).
$$
{\revv Then via the auxiliary function $W$ discussed in the proof of Lemma~\ref{var-minibatch}, we have the following.
\begin{lemma}\label{lm: vrdual1} Consider Algorithm~\ref{alg: DI-SGD} and using variance reduction to estimate $G^k$ \eqref{def: Gk-SPIDER}. \textcolor{black}{Let $p$ and $q$ be H{\"o}lder conjugates and let $W$ be the auxiliary function introduced in the proof of Lemma~\ref{var-minibatch}.} Then we have that for any $k \neq n_k$,
\begin{align}
\label{ineq: vrdual1}
        & \bE[ W^*( B^k(\x^k) - B^k(\x^{k-1}) - \nabla f(\x^k) + \nabla f(\x^{k-1}) ) \mid \x^k, \x^{k-1}] \\
        \notag
        & \le \frac{1}{2m_k}\bE[\| \nabla F(\x^k,\zeta^k_1) - \nabla F(\x^{k-1}, \zeta^k_1) - \nabla f(\x^k) + \nabla f(\x^{k-1}) \|_q^2 \mid \x^k, \x^{k-1}].
\end{align}
\end{lemma}
\begin{proof}
    The proof is similiar to the proof of Lemma~\ref{var-minibatch}. Let us denote for $k \neq n_k$, 
    $$
    \Delta^k_i \triangleq \nabla F(\x^k, \zeta^k_i) - \nabla F(\x^{k-1}, \zeta^k_i) - \nabla f(\x^k) + \nabla f(\x^{k-1}).
    $$
    Then
    \begin{align*}
        B^k(\x^k) - B^k(\x^{k-1}) - \nabla f(\x^k) + \nabla f(\x^{k-1}) = \frac{1}{m_k}\sum_{i=1}^{m_k} \Delta^k_i.
    \end{align*}
    Therefore, by 1-Lipschitz smoothness of $W^*$, we have that
    \begin{align*}
    & W^*( B^k(\x^k) - B^k(\x^{k-1}) - \nabla f(\x^k) + \nabla f(\x^{k-1}) ) \\
    & \le W^*\left( \frac{1}{m_k}\sum_{i=1}^{m_k-1} \Delta^k_i \right) + \left\langle \nabla W^*\left( \frac{1}{m_k}\sum_{i=1}^{m_k-1} \Delta^k_i \right), \frac{1}{m_k} \Delta^k_{m_k} \right\rangle + \frac{1}{2m_k^2} \| \Delta_{m_k}^k \|_q^2.
    \end{align*}
    {\revi To make this iteration explicit, define
    \[
        S_j^k \triangleq \frac{1}{m_k}\sum_{i=1}^j \Delta_i^k,
        \qquad S_0^k \triangleq 0,
    \]
    For each $j=1,\ldots,m_k$, the same smoothness inequality gives
    \begin{align*}
        W^*\left( \frac{1}{m_k}\sum_{i=1}^j \Delta_i^k  \right)
        \le {}& W^*(S_{j-1}^k)
        + \left\langle \nabla W^*(S_{j-1}^k),
        \frac{1}{m_k}\Delta_j^k \right\rangle
        + \frac{1}{2m_k^2}\|\Delta_j^k\|_q^2.
    \end{align*}
    The current samples $\zeta_1^k,\ldots,\zeta_{m_k}^k$ are independent of the algorithmic history and are i.i.d. Hence, conditional on $\x^k,\x^{k-1},\zeta_1^k,\ldots,\zeta_{j-1}^k $, expectation of the inner-product term is zero. Taking conditional expectations, summing from $j=1$ to $m_k$, and using $W^*(S_0^k)=W^*(0)=0$ therefore yield
    \begin{align*}
        \bE\!\left[
        W^*\!\left(\frac{1}{m_k}\sum_{i=1}^{m_k}\Delta_i^k\right)
        \,\middle|\,\x^k,\x^{k-1}\right] & \le
        \frac{1}{2m_k^2}\sum_{j=1}^{m_k}
        \bE\!\left[\|\Delta_j^k\|_q^2
        \,\middle|\,\x^k,\x^{k-1}\right] \\
        & =
        \frac{1}{2m_k}
        \bE\!\left[\|\Delta_1^k\|_q^2
        \,\middle|\,\x^k,\x^{k-1}\right],
    \end{align*}
    where the last equality follows because the $\Delta_j^k$ are conditionally identically distributed. Substituting the definition of $\Delta_1^k$ gives \eqref{ineq: vrdual1}.}
    \qed
\end{proof}

{\revi Before stating the next bound, let $\sF_k$ denote the filtration generated by the algorithmic history up to the beginning of iteration $k$, including $\x^1,\ldots,\x^k$ and all samples used before iteration $k$.  Recall that the random output index $Y$ is the index generated in Step 3 of Algorithm~\ref{alg: DI-SGD}.}
\begin{lemma}\label{lm: spiderbd} 
Consider Algorithm~\ref{alg: DI-SGD}. Suppose that Assumption~\ref{assumption first-order oracle}, \ref{ass: subG} and ~\ref{ass: Lip} hold. Consider using variance reduction to estimate $G^k$ {\revi as in \eqref{def: Gk-SPIDER} and \eqref{mk_spider}.  Let $p_k=1/K$. } Then
\begin{align}
\label{redvar-SPIDER2}
 \bE[ \| \nabla f(\x^Y) - G^Y \|_\infty^2 ]
    \le \frac{4 {\revi e^4 (\log d)^2} L^2 q_0}{m} \bE[ \| \x^{Y+1} - \x^Y \|_1^2 ] + \frac{{\revi e^4 (\log d)^2 }\sigma_\infty^2}{m_1},
\end{align}
\end{lemma}
\begin{proof}
By leveraging the property of $W^*$, we can derive the following for $k \neq n_k$:
\begin{align}
\notag
        & \bE[ W^*(G^k - \nabla f(\x^k)) ] \\
        \notag
        & = \bE[ W^*(B^k(\x^k) - B^k(\x^{k-1}) + G^{k-1} - \nabla f(\x^k) + \nabla f(\x^{k-1}) - \nabla f(\x^{k-1}) ) ] \\
        \notag
        & \le \bE[ W^*( G^{k-1} - \nabla f(\x^{k-1}) ) ] \\
        \notag
        & + \bE[ \left\langle \nabla W^*(G^{k-1} - \nabla f(\x^{k-1})), B^k(\x^k) - B^k(\x^{k-1})  - \nabla f(\x^k) + \nabla f(\x^{k-1}) \right\rangle ]  \\
        \notag
        & + (1/2)\bE[ \| B^k(\x^k) - B^k(\x^{k-1})  - \nabla f(\x^k) + \nabla f(\x^{k-1}) \|_q^2 ] \\
        \notag
        & = \bE[ W^*( G^{k-1} - \nabla f(\x^{k-1}) ) ] \\
        \label{ineq: spiderbd1}
        & + (1/2)\bE[ \| B^k(\x^k) - B^k(\x^{k-1})  - \nabla f(\x^k) + \nabla f(\x^{k-1}) \|_q^2 ],
\end{align}
where the last equality is because of the fact that
\begin{align*}
   &  \bE[ \left\langle \nabla W^*(G^{k-1} - \nabla f(\x^{k-1})), B^k(\x^k) - B^k(\x^{k-1})  - \nabla f(\x^k) + \nabla f(\x^{k-1}) \right\rangle ] \\
   & {\revi = \bE[ \bE[ \left\langle \nabla W^*(G^{k-1} - \nabla f(\x^{k-1})), B^k(\x^k) - B^k(\x^{k-1})  - \nabla f(\x^k) + \nabla f(\x^{k-1}) \right\rangle \mid \mathcal{F}_k ] ] }\\
   & = 0.
\end{align*}
Note that by Lemma~\ref{lm: vrdual1}, we have for any $k \neq n_k$,
\begin{align*}
    & \bE[ \| B^k(\x^k) - B^k(\x^{k-1}) - \nabla f(\x^k) + \nabla f(\x^{k-1}) \|_q^2 ] \\
    & \le 2\Omega \bE[ W^*( B^k(\x^k) - B^k(\x^{k-1}) - \nabla f(\x^k) + \nabla f(\x^{k-1}) ) ] \\
    & \le \frac{\Omega}{m_k} \bE[ \| \nabla F(\x^k,\zeta^k_1) - \nabla F(\x^{k-1},\zeta^{k}_1) - \nabla f(\x^k) + \nabla f(\x^{k-1}) \|_q^2 ].
\end{align*}
By plugging the above into \eqref{ineq: spiderbd1} and iteratively invoking \eqref{ineq: spiderbd1}, we have that
\begin{align*}
\notag
    & \bE[W^*(G^k - \nabla f(\x^k))] \\
    & \le \bE[W^*(G^{n_k} - \nabla f(\x^{n_k}))] \\
    & + \frac{\Omega}{2m} \sum_{i=n_k+1}^k \bE[ \| \nabla F(\x^i,\zeta^i_1) - \nabla F(\x^{i-1},\zeta^i_1) - \nabla f(\x^i) + \nabla f(\x^{i-1}) \|_q^2 ] \\
    & \le \bE[ W^*(0) +  \langle \nabla W^*(0), G^{n_k} - \nabla f(\x^{n_k}) \rangle +  (1/2)  \| G^{n_k} - \nabla f(\x^{n_k}) \|_q^2  \\
    & + \frac{\Omega}{2m} \sum_{i=n_k+1}^k \bE[ \| \nabla F(\x^i,\zeta^i_1) - \nabla F(\x^{i-1},\zeta^i_1) - \nabla f(\x^i) + \nabla f(\x^{i-1}) \|_q^2 ] \\
    & = \frac{1}{2}\bE[ \| G^{n_k} - \nabla f(\x^{n_k}) \|_q^2 ] \\
    & + \frac{\Omega}{2m} \sum_{i=n_k+1}^k \bE[ \| \nabla F(\x^i,\zeta^i_1) - \nabla F(\x^{i-1},\zeta^i_1) - \nabla f(\x^i) + \nabla f(\x^{i-1}) \|_q^2 ].
\end{align*}
Therefore, we have that
\begin{align}
\label{ineq: spiderbd2}
\notag
    & \bE[ \| G^k - \nabla f(\x^k) \|_q^2  ] \\
    \notag 
    & \le 2 \Omega \bE[ W^*(G^k - \nabla f(\x^k)) ] \\
    \notag
& \le \Omega \bE[ \| G^{n_k} - \nabla f(\x^{n_k}) \|_q^2 ] \\
& + \frac{\Omega^2}{m} \sum_{i=n_k+1}^k \bE[ \| \nabla F(\x^i,\zeta^i_1) - \nabla F(\x^{i-1},\zeta^i_1) - \nabla f(\x^i) + \nabla f(\x^{i-1}) \|_q^2 ].
\end{align}

Now let $p \textcolor{black}{\triangleq} 1$ and $q \textcolor{black}{\triangleq} +\infty$. Then by Assumption~\ref{ass: subG}, \ref{ass: Lip}, \eqref{ineq: spiderbd2} and \eqref{redvar-SPIDER1}, we have that
\begin{align}
\notag
    & \bE[ \| G^k - \nabla f(\x^k) \|_{\infty}^2  ] \\
    \notag
& \le \frac{\Omega e^2 (\log d) \sigma_\infty^2}{m_1} \\
\notag
& + \frac{\Omega^2}{m} \sum_{i=n_k+1}^k \bE[ \| \nabla F(\x^i,\zeta^i_1) - \nabla F(\x^{i-1},\zeta^i_1) - \nabla f(\x^i) + \nabla f(\x^{i-1}) \|_{\infty}^2 ] \\
\label{ineq: spiderbd3}
& \le \frac{\Omega e^2 (\log d) \sigma_\infty^2}{m_1} + \frac{4\Omega^2 L^2}{m} \sum_{i= n_k+1}^k \bE[ \| \Delta \x^i \|_1^2 ],
\end{align}
Therefore,
\begin{align}
\notag
    & \bE[ \| \nabla f(\x^Y) - G^Y \|_\infty^2 ] \\
    \notag
    & = \frac{1}{K} \sum_{k=1}^K \bE[ \| \nabla f(\x^k) - G^k \|_\infty^2 ] \\
    \notag
    & \overset{\eqref{ineq: spiderbd3}}{\le} \frac{4 \Omega^2 L^2}{K m} \sum_{k=2,k \neq n_k}^K \sum_{n_k + 1}^k \bE[\| \Delta \x^i \|_1^2] +\frac{\Omega e^2(\log d)\sigma_\infty^2}{m_1} \\
    \notag
    & \le \frac{4\Omega^2 L^2q_0}{Km} \sum_{k=1}^K \bE[\| \x^{k+1} - \x^k \|_1^2] +  \frac{\Omega e^2(\log d)\sigma_\infty^2}{m_1} \\
    \notag
    & = \frac{4\Omega^2 L^2 q_0}{m} \bE[ \| \x^{Y+1} - \x^Y \|_1^2 ] + \frac{\Omega e^2(\log d)\sigma_\infty^2}{m_1}.
\end{align}
Plug in $\Omega = e^2 (\log d)$ and we conclude the result. \qed
\end{proof}
\begin{remark}
    \eqref{redvar-SPIDER2} is vital in complexity analysis when variance reduction \eqref{def: Gk-SPIDER} is applied.  In 2-norm setting,  a similar bound holds, compared to which we have extra scalars in order of $(\log d)^2$ in the bound for the non-Euclidean setting.
\end{remark}}

\section{Convergence and complexity analysis}
\label{sec: analysis}
In this section we analyze Algorithm~\ref{alg: DI-SGD} in the various settings discussed in the last section and derive sample complexity guarantees.  Throughout this section, 
let $\Delta_f \triangleq \bE[f(x^1)] - \underline{f}$ and $\Omega \triangleq e^2 \log d$. 

{\revi Let $X_k \triangleq X + \{ -\x^k \}$ where $+$ is the Minkovski summation. Then by substituting ${\y} \textcolor{black}{\triangleq} \x - \x^k$ and optimality condition for non-smooth convex optimization problem \eqref{def: proj}, we have
\begin{align*}
    0 \in (\hat{\y} + \x^k - \vbf) + \partial (\phi + \delta_{X_k})(\hat{\y}),
\end{align*}
where $\hat{\y} \textcolor{black}{\triangleq} P_X^k(\vbf) - \x^k$.  Note that for our choices of $\phi$, the following condition holds:
\begin{align}\label{ri}
{\rm ri}({\rm dom}\; \phi) \cap {\rm ri} ({\rm dom}\; \delta_{X_k}) \neq \emptyset.
\end{align}
Indeed,  when $\phi(\x) = (\hrho/2)\| \x \|_1^2$,  ${\rm dom} \; \phi = \bR^d$ and \eqref{ri} naturally holds; when $ \phi(\x) = \delta_{\bar X}(\x) $ where $\bar X$ is the $\ell_1$-ball centered at origin with radius $\psi$,  \eqref{ri} holds since $0$ is an interior point of $\bar X$ and $0 \in X_k$, a convex set.  Because \eqref{ri} implies
\begin{align*}
    \partial (\phi + \delta_{X_k})(\cdot) = \partial \phi (\cdot) + \partial \delta_{X_k}(\cdot)
\end{align*}
(see \cite[Thm. 4.5]{drusvyatskiy2020convex}),  we have
\begin{align}\label{optcon-inclu1}
0 \in (\hat{\y} + \x^k - \vbf) + \partial \phi (\hat \y) + \partial \delta_{X_k}(\hat{\y}).
\end{align}
 Note that $\partial \delta_{X_k}(\cdot) = \mathcal{N}_{X_k}(\cdot)$.  Then the following conditions hold by definition of $\x^{k+1}$ and letting $\vbf \textcolor{black}{\triangleq} \x^k - \eta_k G^k$ in \eqref{optcon-inclu1}: $\exists \xi^{k+1}$, s.t.,
\begin{subequations}\label{ineq: proj}
\begin{align}\label{ineq: proj1} 
        (\Delta \xbf^{k+1} + \eta_k G^k + \xi^{k+1})^T( \xbf - \xbf^{k+1} ) \ge 0, &\quad 
        \forall \xbf \in X, \\
        \label{eq: proj2}
        \xi^{k+1} \in \partial \phi(\cdot)(\xbf^{k+1} - \x^k).
    \end{align}
\end{subequations}
In particular, when $X = \bR^d$, we have
\begin{align} \label{eq: proj4}
    - \Delta \x^{k+1} - \eta_k G^k = \xi^{k+1} \in \partial \phi(\cdot)(\x^{k+1} - \x^k).
\end{align}}

Next we provide two key lemmas useful for analyzing all the cases.
\subsection{Key lemmas}
We start with the following lemma which holds without specifying the choices of $\phi$ or $G^k$.
{\revi
\begin{lemma}\label{lm: errbd1} Consider Algorithm~\ref{alg: DI-SGD} with \eqref{ineq: proj}.  Suppose that Assumption~\ref{assumption first-order oracle}, \ref{asp: flip}, \ref{ass: setw3hj} hold. Then for any positive scalar sequences $\{ \eta_k \}_k$ and $\{ t_k \}_k$, we have
\begin{align}
\notag
   & \sum_{k=1}^K  \frac{1}{\eta_k} \bE [ \| \Delta \xbf^{k+1} \|^2 ] - \sum_{k=1}^K \left( \frac{L}{2} + \frac{1}{2 t_k} \right)  \bE [ \| \Delta \xbf^{k+1} \|_1^2 ] + \sum_{k=1}^K \frac{1}{\eta_k} \bE[ ( \xi^{k+1} )^T  \Delta \xbf^{k+1} ] \\
   \label{errbd1}
    & \le \Delta_f + \sum_{k=1}^K  \frac{t_k}{2} \bE[ \| g^k - G^k \|_\infty^2 ].
\end{align}
\end{lemma}
}
\begin{proof}
    We have that 
\begin{align*}
 f(\xbf^{k+1})
    & \overset{\eqref{ineq: Taylor}}{\le} f(\xbf^k) + (g^k)^T \Delta \xbf^{k+1} + \frac{L}{2} \| \Delta \xbf^{k+1} \|_1^2 \\
    & \le f(\xbf^k) + (g^k - G^k)^T \Delta \xbf^{k+1} + (G^k)^T \Delta \xbf^{k+1} +  \frac{L}{2} \| \Delta \xbf^{k+1} \|_1^2 \\
    & \overset{\eqref{ineq: proj1}}{\le} f(\xbf^k) + (g^k - G^k)^T \Delta \xbf^{k+1} \\
    & + \frac{1}{\eta_k}(\Delta \xbf^{k+1} + \xi^{k+1} )^T(\xbf^k - \xbf^{k+1}) +  \frac{L}{2} \| \Delta \xbf^{k+1} \|_1^2 \\
    & \le f(\xbf^k) + \frac{t_k}{2} \| g^k - G^k \|_\infty^2 + \frac{1}{2 t_k} \| \Delta \xbf^{k+1} \|_1^2 - \frac{1}{\eta_k} \| \Delta \xbf^{k+1} \|^2  \\
    & + \frac{1}{\eta_k} ( \xi^{k+1} )^T  (\xbf^k - \xbf^{k+1}) + \frac{L}{2}  \| \Delta \xbf^{k+1} \|_1^2.
\end{align*}
By rearranging terms, we have that
\begin{align*}
    & \frac{1}{\eta_k} \| \Delta \xbf^{k+1} \|^2 - \left( \frac{L}{2} + \frac{1}{2 t_k} \right) \| \Delta \xbf^{k+1} \|_1^2 + \frac{1}{\eta_k} ( \xi^{k+1} )^T  (\xbf^{k+1} - \xbf^k) \\
    & \le f(\xbf^k) - f(\xbf^{k+1}) + \frac{t_k}{2} \| g^k - G^k \|_\infty^2.
\end{align*}
Sum up the above inequality for $k = 1,\hdots,K$ and take expectation on both sides, we have
\begin{align*}
   & \sum_{k=1}^K  \frac{1}{\eta_k} \bE [ \| \Delta \xbf^{k+1} \|^2 ] - \sum_{k=1}^K \left( \frac{L}{2} + \frac{1}{2 t_k} \right)  \bE [ \| \Delta \xbf^{k+1} \|_1^2 ] + \sum_{k=1}^K \frac{1}{\eta_k} \bE[ ( \xi^{k+1} )^T  \Delta \xbf^{k+1} ] \\
    & \le \bE[ f(\xbf^1) - f(\xbf^{K+1}) ] + \sum_{k=1}^K  \frac{ t_k}{2} \bE[ \| g^k - G^k \|_\infty^2 ].
\end{align*}
Recall $\Delta_f = \bE[f(\x^1)] - \lf \ge \bE[ f(\xbf^1) - f(\xbf^{K+1}) ]$. Then the result follows.\qed
\end{proof}

The residual of interest is ${\rm dist}_{\| \cdot \|_\infty}(0, \partial(f+\delta_X)(\x^{Y+1})) $. The following lemma provides an upper bound for its expected value.
{\revi
\begin{lemma} Consider Algorithm~\ref{alg: DI-SGD} with \eqref{ineq: proj}.  Suppose that Assumption~\ref{assumption first-order oracle}, \ref{asp: flip} hold.  For any positive stepsize sequence $\{\eta_k\}_k$,  we have
    \begin{align}
\label{errbd2}
    & \bE[ {\rm dist}_{\| \cdot \|_\infty}(0, \partial(f+\delta_X)(\x^{Y+1}))] \\
\notag
    & \le L \bE[ \| \Delta \x^{Y+1} \|_1 ] + \bE[ \| g^Y - G^Y \|_\infty ] + \bE[ \| \Delta \x^{Y+1} \|_\infty / \eta_Y ] + \bE[ \| \xi^{Y+1} \|_\infty / \eta_Y ].
\end{align}
\end{lemma}
}
\begin{proof}
Note that due to \eqref{ineq: proj1}, we have that
\begin{align*}
& - (\Delta \x^{Y+1} + \eta_k G^Y + \xi^{Y+1})  \in N_X(\x^{Y+1}) \\
\implies & - \frac{1}{\eta_k}(\Delta \x^{Y+1} + \eta_k G^Y + \xi^{Y+1} )  \in \partial \delta_X(\x^{Y+1}) \\
\implies &  \nabla f(\x^{Y+1}) - \frac{1}{\eta_k}(\Delta \x^{Y+1} + \eta_k G^Y + \xi^{Y+1} ) \in \partial(f+\delta_X)(\x^{Y+1}),
\end{align*}
where the last implication is due to differentiability of $f$. 
Therefore, 
\begin{align}
\notag
    & \bE[ {\rm dist}_{\| \cdot \|_\infty}(0, \partial(f+\delta_X)(\x^{Y+1}))] \\
    \notag
    & \le \bE[\| \nabla f(\x^{Y+1}) - \frac{1}{\eta_Y}(\Delta \x^{Y+1} + \eta_Y G^Y + \xi^{Y+1} ) \|_\infty] \\
    \notag
    & = \bE[ \|  \nabla f(\x^{Y+1}) - g^Y + g^Y - G^Y - (\Delta \x^{Y+1}  + \xi^{Y+1} )/\eta_Y  \|_\infty ] \\
    \notag
    & \le \bE[ \| \nabla f(\x^{Y+1}) - g^Y \|_\infty ]  + \bE[ \| g^Y - G^Y \|_\infty] + \bE\left[ \| \Delta \x^{Y+1} + \xi^{Y+1} \|_\infty / \eta_Y \right] \\\notag
    & \overset{\revv \eqref{ineq: lip}}{\le} L \bE[ \| \Delta \x^{Y+1} \|_1 ] + \bE[ \| g^Y - G^Y \|_\infty ] + \bE[ \| \Delta \x^{Y+1} \|_\infty / \eta_Y ] + \bE[ \| \xi^{Y+1} \|_\infty / \eta_Y ] .  
\end{align}\qed
\end{proof}
\subsection{Analysis for Case 1}
The main result for minibatch sampling is as follows.
{\revi
\begin{theorem}\label{thm: case1mini} For Algorithm~\ref{alg: DI-SGD} with \eqref{ineq: proj}, consider case 1 \eqref{setting1} and estimating $G^k$ using minibatch \eqref{def: Gk-SGD}. {Suppose that Assumption~\ref{assumption first-order oracle}, \ref{ass: subG}, \ref{asp: flip} and \ref{ass: setw3hj} hold}.  Suppose that $\hrho > 3/2$,  $ 0 < \eta_k \le \frac{1}{L} $, $t_k \triangleq \frac{3 \eta_k}{\hrho}$  and let $p_k \triangleq \frac{\eta_k}{\sum_{k=1}^K \eta_k}$.
Then
\begin{align}
\notag
    & \bE[ {\rm dist}_{\| \cdot \|_\infty}(0, \partial(f+\delta_X)(\x^{Y+1})) ] \\
    & \le \sqrt{ \frac{ c_1(\hrho) \Delta_f }{\sum_{k=1}^K \eta_k} } +  \sqrt{ \frac{ \left( \frac{ 3 c_1(\hrho)}{ 2 \hrho} +1 \right) \sum_{k=1}^K \frac{\eta_k}{ m_k}  }{ \sum_{k=1}^K \eta_k   } } \cdot \sqrt{e^2  (\log d) \sigma_\infty^2 },
    \label{errbd4}
\end{align}
where $c_1(\hrho)$ is a constant depending on $\hrho$.
\end{theorem}}
\begin{proof}
By Lemma~\ref{lm: l1subdiff}
we have that
\begin{align*}
    \partial \phi(\z - \x^k) \subseteq \{ \hrho \| \z - \x^k \|_1 v \mid v_i \in [-1,1] \}.
\end{align*}
Therefore,  by \eqref{eq: proj2},  
\begin{align}\label{subgbd}
    \| \xi^{k+1} \|_\infty \le \hrho \| \x^{k+1} - \x^k \|_1.
\end{align}
Note that $\phi(0) = 0$, then by \eqref{eq: proj2} and convexity of $\phi$,
\begin{align}\label{ineq: phi}
    ( \xi^{k+1} )^T  (\xbf^k - \xbf^{k+1}) + \phi(\xbf^{k+1}-\x^k) \le \phi(0) = 0.
\end{align}

Since we consider using a minibatch of samples to estimate the gradient \eqref{def: Gk-SGD}. Then \eqref{errbd1}\eqref{setting1}\eqref{subgbd}\eqref{ineq: phi} lead to
\begin{align}
\notag
& \sum_{k=1}^K \bE \left[ \frac{1}{\eta_k} \| \Delta \xbf^{k+1} \|^2 \right] - \sum_{k=1}^K \left( \frac{L}{2} + \frac{1}{2 t_k} \right) \bE [ \| \Delta \xbf^{k+1} \|_1^2 ] + \sum_{k=1}^K \frac{\hrho}{2\eta_k} \bE[ \| \Delta \xbf^{k+1} \|_1^2 ] \\\notag
& \le \sum_{k=1}^K \bE \left[ \frac{1}{\eta_k} \| \Delta \xbf^{k+1} \|^2 \right] - \sum_{k=1}^K \left( \frac{L}{2} + \frac{1}{2 t_k} \right) \bE [ \| \Delta \xbf^{k+1} \|_1^2 ] \\ \notag
& + \sum_{k=1}^K \frac{1}{\eta_k} \bE[ ( \xi^{k+1} )^T \Delta \x^{k+1}  ] \\\notag
& \le \Delta_f + \sum_{k=1}^K \frac{t_k}{2} \bE[ \| g^k - G^k \|_\infty^2 ] \\
\notag
& \implies  \sum_{k=1}^K \frac{1}{\eta_k} \bE [ \| \Delta \xbf^{k+1} \|^2 ] +  \sum_{k=1}^K \left( \frac{\hat \rho}{2 \eta_k } - \frac{L}{2} -  \frac{1}{2 t_k}  \right) \bE[ \| \Delta \xbf^{k+1} \|_1^2 ] \\
\label{errbd3}
& \le \Delta_f +  \sum_{k=1}^K \frac{t_k}{2} \bE[ \| g^k - G^k \|_\infty^2 ]. 
\end{align}
{\revi
Let us denote $\alpha_k \triangleq \frac{ \frac{\hat{\rho}}{2\eta_k} - \frac{L}{2} - \frac{1}{2t_k} }{ \left( L + (\hat{\rho}+1)/\eta_k \right)^2} $,  then we have
\begin{align*}
\alpha_k 
&= \frac{ \frac{\hat{\rho}}{3\eta_k} - \frac{L}{2} }{ \left( L + (\hat{\rho}+1)/\eta_k \right)^2} = \frac{ \frac{\hat{\rho}}{3} \eta_k - \frac{L}{2} \eta_k^2 }{ (L\eta_k + \hat{\rho} + 1)^2 }.
\end{align*}
Since $\eta_k \le 1/L$,  if we denote $c_1(\hrho) \triangleq (\hrho + 2)^2/(\hrho/3 - 1/2)$, then
\begin{align}\label{ineq: alphak}
\alpha_k \geq \eta_k \left( \hat{\rho}/3 - 1/2 \right) \big/ (\hat{\rho} + 2)^2 = \eta_k/c_1(\hat{\rho}).
\end{align}
By \eqref{errbd2}\eqref{subgbd}\eqref{errbd3}\eqref{ineq: alphak} and Lemma~\ref{var-minibatch},
we have
\begin{align*}
    & \bE[ {\rm dist}_{\| \cdot \|_\infty}(0, \partial(f+\delta_X)(\x^{Y+1}))] \\
\notag
 & \le L \bE[ \| \Delta \x^{Y+1} \|_1 ] + \bE[ \| g^Y - G^Y \|_\infty ] + \bE[ \| \Delta \x^{Y+1} \|_\infty / \eta_Y ] + \bE[ \| \xi^{Y+1} \|_\infty / \eta_Y ] \\
 & \le  \bE\left[ \left( L +  \frac{\hat \rho + 1}{\eta_Y} \right) \| \Delta \x^{Y+1} \|_1 \right] + \bE[ \| g^Y - G^Y \|_\infty ] \\
 & \le \sqrt{ \bE\left[ \left( L +  \frac{\hat \rho + 1}{\eta_Y} \right)^2 \| \Delta \x^{Y+1} \|_1^2 \right] } + \bE[ \| g^Y - G^Y \|_\infty ] \\
 & =  \sqrt{ \sum_{k=1}^K  \bE\left[ \left( L +  \frac{\hat \rho + 1}{\eta_k} \right)^2 \| \Delta \x^{k+1} \|_1^2 \right] P(Y = k) } + \bE[ \| g^Y - G^Y \|_\infty ] \\
   & =  \sqrt{ \frac{ \sum_{k=1}^K \left( L +  \frac{\hat \rho + 1}{\eta_k} \right)^2 \eta_k \bE\left[ \| \Delta \x^{k+1} \|_1^2 \right] }{\sum_{k=1}^K \eta_k} } + \bE[ \| g^Y - G^Y \|_\infty ] \\
  & \le \sqrt{ \frac{  c_1(\hrho) \sum_{k=1}^K \left( \frac{\hat \rho}{2 \eta_k} - \frac{L}{2} - \frac{1}{2 t_k} \right)  \bE\left[ \| \Delta \x^{k+1} \|_1^2 \right] }{\sum_{k=1}^K \eta_k} } + \bE[ \| g^Y - G^Y \|_\infty ] \\
  & \le  \sqrt{ \frac{ c_1(\hrho) \left( \Delta_f + \frac{1}{2} \sum_{k=1}^K t_k \bE[\| g^k - G^k \|_\infty^2 ] \right) }{\sum_{k=1}^K \eta_k} } + \bE[ \| g^Y - G^Y \|_\infty ] \\
  & \le  \sqrt{ \frac{ c_1(\hrho) \Delta_f }{\sum_{k=1}^K \eta_k} } + \sqrt{ \frac{ \frac{c_1(\hrho)}{2} \sum_{k=1}^K t_k \bE[\| g^k - G^k \|_\infty^2 ] }{\sum_{k=1}^K \eta_k} } + \bE[ \| g^Y - G^Y \|_\infty ] \\
  &  \le \sqrt{ \frac{ c_1(\hrho) \Delta_f }{\sum_{k=1}^K \eta_k} } + \sqrt{ \frac{ \frac{c_1(\hrho)}{2} \sum_{k=1}^K t_k \bE[\| g^k - G^k \|_\infty^2 ] }{\sum_{k=1}^K \eta_k} } + \sqrt{ \frac{ \sum_{k=1}^K \eta_k \bE[\| g^k - G^k \|_\infty^2 ] }{\sum_{k=1}^K \eta_k} } \\
  &  \le \sqrt{ \frac{ c_1(\hrho) \Delta_f }{\sum_{k=1}^K \eta_k} } +  \sqrt{ \frac{ \left( \frac{3c_1(\hrho)}{2 \hrho} +1 \right) \sum_{k=1}^K \frac{\eta_k}{ m_k}  }{ \sum_{k=1}^K \eta_k   } } \cdot \sqrt{e^2  (\log d) \sigma_\infty^2 } .
\end{align*}
}
\qed
\end{proof}
\begin{corollary}\label{corr: mb.comp.case1}
Under the setting of Theorem~\ref{thm: case1mini}, let
\begin{align*}
    \eta_k \equiv \frac{1}{L},\; K \triangleq \left\lceil \frac{\Delta_f L}{\epsilon^2} \right\rceil,\; m_k \equiv \left\lceil \frac{e^2 (\log d) \sigma_\infty^2}{\epsilon^2} \right\rceil.  
\end{align*}
Then we have
\begin{align*}
\bE[ {\rm dist}_{\| \cdot \|_\infty}(0, \partial(f+\delta_X)(\x^{Y+1})) ] \le {\revv C_1(\hrho) \epsilon},
\end{align*}
where $C_1(\hat \rho)$ is a constant depending only on $\hat \rho$.  As a result,  the sample complexity of the algorithm is
\begin{align*}
    Km = \left\lceil \frac{\Delta_f L}{\epsilon^2} \right\rceil \left\lceil \frac{e^2 (\log d) \sigma_\infty^2}{\epsilon^2} \right\rceil = \sO \left(\frac{\Delta_f L (\log d) \sigma_\infty^2}{\epsilon^4} \right).
\end{align*}
\end{corollary}
{\revi
To facilitate diminishing stepsize and asymptotic convergence, we consider the following setting.
\begin{corollary}\label{corr: case1.asymp}
Under the setting of Theorem \ref{thm: case1mini}, let $\alpha \in (0,1)$ and 
\begin{align*}
\eta_k = k^{-\alpha} / L,  \; m_k = \lceil m_1 k \rceil.
\end{align*}
Then we have
\begin{align*}
& \mathbb{E}[ {\rm dist}_{ \| \cdot \|_\infty} (0, \partial(f+\delta_X)(x^{Y+1})) ] \\
& \le \sqrt{ \frac{ c_1(\hat{\rho}) \Delta_f L(1-\alpha) }{(K+1)^{1-\alpha} - 1} } + \sqrt{ \frac{ \left( \frac{ 3 c_1(\hat{\rho})}{ 2 \hat{\rho}}+1 \right)(1/\alpha - \alpha)}{(K+1)^{1-\alpha} - 1} } \cdot \sqrt{ \frac{ e^2 (\log d) \sigma_\infty^2}{m_1} },
\end{align*}
where $c_1(\hat{\rho})$ is a constant depending only on $\hat{\rho}$.
\end{corollary}
\begin{proof}
By plugging the parameter selections into \eqref{errbd4}, we have
\begin{align*}
&\mathbb{E}[ \text{dist}_{\|\cdot\|_\infty} (0, \partial(f+\delta_X)(\x^{Y+1})) ] \\
&\le \sqrt{ \frac{c_1(\hat{\rho}) \Delta_f L}{\sum_{k=1}^K k^{-\alpha}} } + \sqrt{ \frac{ \left( \frac{3c_1(\hat{\rho})}{2\hat{\rho}} + 1 \right) \sum_{k=1}^K k^{-(\alpha+1)}}{\sum_{k=1}^K k^{-\alpha}} } \cdot \sqrt{\frac{e^2 (\log d) \sigma_\infty^2}{m_1}} \\
&\le \sqrt{ \frac{ c_1(\hat{\rho}) \Delta_f L(1-\alpha) }{(K+1)^{1-\alpha} - 1} } + \sqrt{ \frac{ \left( \frac{3c_1(\hat{\rho})}{2\hat{\rho}}+1 \right)(1/\alpha - \alpha)}{(K+1)^{1-\alpha} - 1} } \cdot \sqrt{ \frac{e^2 (\log d) \sigma_\infty^2}{m_1} }.
\end{align*}
where the second inequality leverages the facts that
\begin{align*}
\sum_{k=1}^K k^{-\alpha} &\ge \int_1^{K+1} x^{-\alpha} dx = \frac{x^{1-\alpha}}{1-\alpha} \bigg|_1^{K+1} = \frac{(K+1)^{1-\alpha} - 1}{1-\alpha}, \\
\sum_{k=1}^K k^{-(\alpha+1)} &\le \int_1^K x^{-(\alpha+1)} dx + 1 = \frac{x^{-\alpha}}{-\alpha} \bigg|_1^K + 1 = -\frac{1}{\alpha} \left( \frac{1}{K^\alpha} - 1 \right) + 1 \le 1 + \frac{1}{\alpha}.
\end{align*}
\qed
\end{proof}

\begin{remark} (i). Under the setting of Corollary~\ref{corr: case1.asymp},  the following asymptotic convergence result holds,
\[
\lim_{K \to +\infty} \mathbb{E} [ \text{dist}_{ \| \cdot \|_\infty} (0, \partial(f+\delta_X)(\x^{Y+1})) ]  = 0.
\]
Therefore, 
\[
\liminf_{k \to +\infty} \mathbb{E} [ \text{dist}_{ \| \cdot \|_\infty} (0, \partial(f+\delta_X)(\x^{k+1})) ]  = 0
\]
since
\[
\min\limits_{1 \le k \le K}  \mathbb{E} [ \text{dist}_{ \| \cdot \|_\infty} (0, \partial(f+\delta_X)(\x^{k+1})) ] \le \mathbb{E} [ \text{dist}_{ \| \cdot \|_\infty} (0, \partial(f+\delta_X)(\x^{Y+1})) ].
\]
(ii).  Under the setting of Corollary~\ref{corr: case1.asymp} and $\epsilon  \in (0,1)$,  the algorithm will output a point such that 
$$ \mathbb{E} [ \text{dist}_{ \| \cdot \|_\infty} (0, \partial(f+\delta_X)(\x^{Y+1})) ] \le C_2(\hrho) \epsilon $$ where $C_2(\hrho)$ is a constant depending on the choice of $\hrho$ if
$$
K \ge \frac{ \left[ \max \left\{ \Delta_f L (1-\alpha), (1/\alpha - \alpha) e^2 (\log d) \sigma_\infty^2/ m_1 \right\} + 1 \right]^{\frac{1}{1-\alpha}}}{\epsilon^{\frac{2}{1-\alpha}}} - 1.
$$
Then the sample complexity guarantee under this setting is
\begin{align*}
\sum_{k=1}^K m_k = \sO( m_1 K^2 ) = \tilde{ \sO} \left( \epsilon^{-\frac{4}{1-\alpha}} \right),
\end{align*}
where $\tilde{O}$ hides the polynomial of the $\log d$ term.
\qed
\end{remark}
}
The next theorem reveals the convergence rate by applying variance reduction in Case 1. The proof is similar to Theorem~\ref{thm: case1mini} except for using Lemma~\ref{lm: spiderbd} for complexity improvement and relegated to Appendix~\ref{app: proof}.
\begin{theorem}\label{thm: case1svrg}
For Algorithm~\ref{alg: DI-SGD} with \eqref{ineq: proj}, consider case 1 \eqref{setting1} and estimating $G^k$ using variance reduction {\revv \textcolor{black}{as in \eqref{def: Gk-SPIDER} and \eqref{mk_spider}}}. {\revv Suppose that Assumption~\ref{assumption first-order oracle},~\ref{ass: subG},~\ref{ass: Lip},~\ref{ass: setw3hj} hold.} Suppose that $\eta_k \equiv \eta > 0$,  $ t \triangleq 2\eta/\hrho$,  $\hrho/(4\eta) - L/2  - \frac{2t \Omega^2 L^2 q_0}{m} > 0$ and $p_k \triangleq 1/K$. Then the following holds.
\begin{align}
\notag
    & \bE[ {\rm dist}_{\| \cdot \|_\infty}(0, \partial(f+\delta_X)(\x^{Y+1})) ] \\
    & \le \left( L + \frac{\hrho + 1}{\eta} + \sqrt{\frac{4\Omega^2 L^2q_0}{m}} \right) \sqrt{\frac{\frac{\Delta_f}{K} + \frac{ t \Omega e^2 (\log d)\sigmainf^2}{2m_1} }{\frac{\hrho}{4\eta} - \frac{L}{2} - \frac{2 t\Omega^2L^2q_0}{m} }} + \sqrt{ \frac{\Omega e^2 (\log d)\sigma_\infty^2}{m_1} }.
\label{errbd8}
\end{align}
\end{theorem}
\begin{corollary}\label{corr: vr.comp.case1}
Under the setting of Theorem~\ref{thm: case1svrg}, let
\begin{align*}
    & \eta \triangleq \frac{1}{L},\; K \triangleq \left\lceil \frac{\Delta_f L}{\epsilon^2} \right\rceil,\; m_1 \triangleq \left\lceil \frac{ \Omega e^2 (\log d) \sigma_\infty^2}{\epsilon^2} \right\rceil,\; m \triangleq q_0 \Omega^2,  \; mq_0 = m_1,\; \hrho \ge 6.
\end{align*}
Then by \eqref{errbd8} we have
\begin{align*}
\bE[ {\rm dist}_{\| \cdot \|_\infty}(0, \partial(f+\delta_X)(\x^{Y+1})) ] \le C_3(\hrho)\epsilon,
\end{align*}
where $C_3(\hrho)$ is a constant only dependent on $\hrho$. The sample complexity is upper bounded by
{\revv
\begin{align*}
   Km + \left\lceil \frac{K}{q} \right\rceil m_1 & \le  3Km = \sO\left( \frac{ \Delta_f L \Omega^2 \sigma_\infty }{\epsilon^3} \right) {\revi = \sO\left( \frac{ \Delta_f L (\log d)^2 \sigma_\infty }{\epsilon^3} \right)}.
\end{align*}
}
\end{corollary}

{\revi
\begin{remark}\label{rmk: Lip.const.}
Now let us temporarily denote the Lipschitz constant in \eqref{ineq: lip} as $L_{\infty,1}$ and the one in the Euclidean setting as $L_{2,2}$.  Note that a basic fact is that:
\begin{align*}
\| \nabla f(x) - \nabla f(y) \|_2 \le L_{2,2} \| x- y \|_2 \implies  \| \nabla f(x) - \nabla f(y) \|_\infty \le L_{2,2} \| x- y \|_1.
\end{align*}
Therefore, if both $L_{\infty,1}$ and $L_{2,2}$ are tight estimates of the Lipschitz smoothness, then we have $L_{\infty, 1} \le L_{2,2}$.  Since the stepsize is usually proportional to the reciprocal of the Lipschitz smoothness constant,  smaller $L$ leads to a larger feasible stepsize, thus less tuning efforts.  If $L_{\infty,1}$ is significantly less than $L_{2,2}$ (for e.g., less dimension-sensitive),  then we potentially have better sample complexity,  because the sample complexity's dependence on the Lipschitz constant is linear in either the Euclidean or the non-Euclidean setting.  In fact,  for some problem classes,  $L_{\infty,1}$ can be much less dimension-sensitive than $L_{2,2}$ (see Appendix~\ref{app: Lip}).  
Even when $L_{\infty,1}$ and $L_{2,2}$ are of similar magnitude, we still expect the proposed method to be better in practice in terms of controlling convergence.  Note that stochastic-gradient variance is measured in $\ell_\infty$ for the proposed method, and under componentwise sub-Gaussian noise it is typically logarithmic in $d$, whereas the Euclidean variance is typically linear in $d$.  The former requires much less samples to facilitate convergence when we solve high dimensional problems and the componentwise sub-Gaussian models are appropriate.  
\end{remark}}

\subsection{Analysis for Case 2}
The main result is given as follows.
{\revi
\begin{theorem}\label{thm: case2mini} For Algorithm~\ref{alg: DI-SGD} with \eqref{ineq: proj}, consider {case 2 \eqref{case2phi}} and estimating $G^k$ using minibatch \eqref{def: Gk-SGD}. {Suppose that Assumption~\ref{assumption first-order oracle},~\ref{ass: subG},~\ref{asp: flip},~\ref{ass: setw3hj} hold.} Let $\eta_k \equiv \eta > 0$.  Replace $\psi$ with $\psi_k > 0$, which is allowed to change across iterations and let $p_k \triangleq \frac{\psi_k}{\sum_{k=1}^K \psi_k}$.  Then for any $t > 0$, we have
\begin{align}
\notag
    & \bE[ {\rm dist}_{\| \cdot \|_\infty}(0, \partial(f+\delta_X)(\x^{Y+1}))] \\
    & \le \left(\frac{3L}{2} + \frac{1}{\eta} \right) \frac{\sum_{k=1}^K \psi_k^2}{\sum_{k=1}^K \psi_k} + \frac{\Delta f}{\sum_{k=1}^K \psi_k} + \frac{ 2 \sum_{k=1}^K \frac{ \psi_k }{ \sqrt{m_k} }}{\sum_{k=1}^K \psi_k} \sqrt{e^2 (\log d) \sigma_\infty^2}.
     \label{errbd6}
\end{align}
\end{theorem}
}
\begin{proof}
Note that in the theorem settings, \eqref{eq: proj2} is equivalent to
\begin{align}\label{optcon: subprob}
    \x^{k+1} \in \mbox{arg}\min_\x -(\xi^{k+1})^T \x \quad \mbox{s.t. } \quad \| \x - \x^k \|_1 \le \psi_k.
\end{align}
Therefore, by invoking KKT conditions of \eqref{optcon: subprob}, there exists $\hrho_k$ such that
\begin{align}\label{setting2}
\xi^{k+1} \in \hrho_k \partial (\| \cdot - \x^k \|_1)(\x^{k+1}), \quad 0 \le \hrho_k \perp \| \x^{k+1} - \x^k \|_1 - \psi_k \le 0.
\end{align}
Then
\begin{align*}
    \xi^{k+1} \in \{ \hrho_k v \mid v_i \in [-1,1] \}.
\end{align*}
Therefore, 
\begin{align}\label{subgbd2}
    \| \xi^{k+1} \|_\infty \le \hrho_k.
\end{align}
Moreover, by \eqref{setting2},
\begin{align}\label{ineq: nn3}
    (\xi^{k+1})^T(\x^{k+1} - \x^k) \ge \hrho_k \| \x^{k+1} - \x^k \|_1 - \hrho_k \| \x^k - \x^k \|_1 = \hrho_k \| \x^{k+1} - \x^k \|_1.
\end{align}
{\revi
Consider the minibatch sampling \eqref{def: Gk-SGD}. Under the parameter setting that $\eta_k \equiv \eta$,  \eqref{errbd1}\eqref{setting2}\eqref{subgbd2}\eqref{ineq: nn3} lead to
\begin{align}
\notag
& \sum_{k=1}^K \frac{1}{\eta} \mathbb{E}[\|\Delta \x^{k+1}\|^2] - \sum_{k=1}^K \left(\frac{L}{2} + \frac{1}{2t_k}\right) \mathbb{E}[\|\Delta \x^{k+1}\|_1^2] + \sum_{k=1}^K \mathbb{E}\left[\frac{1}{\eta} \hat{\rho}_k \|\Delta \x^{k+1}\|_1\right] \\
\notag
&  \le{} \sum_{k=1}^K \frac{1}{\eta} \mathbb{E}[\|\Delta \x^{k+1}\|^2] - \sum_{k=1}^K \left(\frac{L}{2} + \frac{1}{2t_k}\right) \mathbb{E}[\|\Delta \x^{k+1}\|_1^2] + \sum_{k=1}^K \frac{1}{\eta} \mathbb{E}\left[(\xi^{k+1})^T \Delta \x^{k+1}\right] \\
\notag
&  \le{} \Delta_f + \sum_{k=1}^K \frac{t_k}{2} \mathbb{E}[\|g^k - G^k\|_\infty^2] \\
\notag
 & \implies  \frac{1}{\eta} \mathbb{E}[\hat{\rho}_Y] = \frac{1}{\eta} \frac{\sum_{k=1}^K \mathbb{E}[\hat{\rho}_k] \cdot \psi_k}{\sum_{k=1}^K \psi_k} \\
 \notag
 & \overset{(Lemma~\ref{var-minibatch})}{\le}  \frac{\Delta_f}{\sum_{k=1}^K \psi_k} + \frac{ L \sum_{k=1}^K \psi_k^2}{ 2 \sum_{k=1}^K \psi_k} +  \frac{\sum_{k=1}^K \left( t_k \cdot \frac{e^2 (\log d) \sigma_\infty^2}{m_k} + \psi_k^2/t_k \right) }{ 2 \sum_{k=1}^K \psi_k}  \\
  \label{errbd5}
 & = \frac{\Delta_f}{\sum_{k=1}^K \psi_k} + \frac{ L \sum_{k=1}^K \psi_k^2}{ 2 \sum_{k=1}^K \psi_k} + \frac{ \sum_{k=1}^K \frac{ \psi_k }{ \sqrt{m_k}}}{ \sum_{k=1}^K \psi_k} \sqrt{e^2 (\log d) \sigma_\infty^2},
\end{align}
where we let $t_k \triangleq \frac{\psi_k \sqrt{m_k}}{ \sqrt{ e^2 (\log d) \sigma_\infty^2} }$ to obtain the last inequality. Therefore, by \eqref{errbd2}\eqref{subgbd2}\eqref{errbd5} and Lemma~\ref{var-minibatch},  we have the following,
\begin{align*}
& \bE[ \text{dist}_{\|\cdot\|_\infty} (0, \partial(f + \delta_X)(\x^{Y+1})) ] \\
\le{} & (L + 1/\eta) \mathbb{E}[\|\Delta \x^{Y+1}\|_1] + \mathbb{E}[\|g^Y - G^Y\|_\infty] + ( 1/\eta ) \mathbb{E}[\|\xi^{Y+1}\|_\infty] \\
\le{} & (L + 1/\eta) \mathbb{E}[\psi_Y] + \frac{\sum_{k=1}^K \mathbb{E}[\|g^k - G^k\|_\infty] \psi_k}{\sum_{k=1}^K \psi_k}  + \frac{1}{\eta} \mathbb{E}[\hat{\rho}_Y] \\
\le{} & (L + 1/\eta) \frac{\sum_{k=1}^K \psi_k^2}{\sum_{k=1}^K \psi_k} + \frac{\sum_{k=1}^K \sqrt{ \mathbb{E}[\|g^k - G^k\|_\infty^2 ] } \psi_k}{\sum_{k=1}^K \psi_k} \\
& +  \frac{\Delta_f}{\sum_{k=1}^K \psi_k} + \frac{ L \sum_{k=1}^K \psi_k^2}{ 2 \sum_{k=1}^K \psi_k} + \frac{ \sum_{k=1}^K \frac{ \psi_k }{ \sqrt{m_k}}}{ \sum_{k=1}^K \psi_k} \sqrt{e^2 (\log d) \sigma_\infty^2} \\
={} & \left(\frac{3L}{2} + \frac{1}{\eta} \right) \frac{\sum_{k=1}^K \psi_k^2}{\sum_{k=1}^K \psi_k} + \frac{\Delta_f}{\sum_{k=1}^K \psi_k} + \frac{ 2 \sum_{k=1}^K \frac{ \psi_k }{ \sqrt{m_k} }}{\sum_{k=1}^K \psi_k} \sqrt{e^2 (\log d) \sigma_\infty^2}.
\end{align*}
 \qed}
\end{proof}
Based on Theorem~\ref{thm: case2mini},  we immediately have the following sample complexity guarantees under an appropriate parameter setting. 
\begin{corollary}\label{corr: mb.comp.case2}
    Under the setting of Theorem~\ref{thm: case2mini}, let
\begin{align*}
    & \eta \triangleq \frac{1}{L},\; K \triangleq  \left\lceil \frac{\Delta_f L}{\epsilon^2} \right\rceil,\; m_k \equiv \left\lceil \frac{e^2 (\log d) \sigma_\infty^2}{\epsilon^2} \right\rceil,\; \psi_k \equiv \frac{\epsilon}{L}.
\end{align*}
Then we have
\begin{align*}
\bE[ {\rm dist}_{\| \cdot \|_\infty}(0, \partial(f+\delta_X)(\x^{Y+1})) ] \le C_1 \epsilon,
\end{align*}
where $C_1$ is an absolute constant. Therefore, the sample complexity is
\begin{align*}
    Km = \sO\left( \frac{ \Delta_f L (\log d) \sigma_\infty^2 }{\epsilon^4} \right).
\end{align*}
\end{corollary}
{\revi
To enable diminishing trust region and asymptotic convergence, we consider the following setting.
\begin{corollary}\label{corr: case2.asymp}
Under the setting of Theorem \ref{thm: case2mini},  suppose that 
\begin{align}\label{param.set}
\psi_k \triangleq  k^{-1/2} / L,  \; m_k \triangleq  \lceil m_1 k \rceil, \; \eta \triangleq  \max\{ 1/L, 1 \}.
\end{align}
Then we have
\begin{align*}
& \mathbb{E}[ {\rm dist}_{ \| \cdot \|_\infty} (0, \partial(f+\delta_X)(\x^{Y+1})) ] \\
& \le \left( \frac{5}{4} + \sqrt{\frac{ e^2 (\log d)\sigma_\infty^2 }{ m_1 }}  \right) \frac{ \log K + 1 }{ \sqrt{K+1} - 1 } + \frac{\Delta_f L}{ 2 ( \sqrt{K+1} - 1) }.
\end{align*}
\end{corollary}
\begin{proof}
This is by plugging the parameter setting \eqref{param.set} in \eqref{errbd6} and the facts that
\begin{align*}
\sum_{k=1}^K \frac{1}{k} & \le 1 + \int_1^K \frac{1}{x} dx = 1 + \log x \bigg|_1^K = \log K + 1, \\
\sum_{k=1}^K \frac{1}{\sqrt{k}} & \ge \int_1^{K+1} \frac{1}{\sqrt{x}} dx = 2 \sqrt{x} \bigg|_1^{K+1} = 2 (\sqrt{K+1} - 1).
\end{align*}\qed
\end{proof}
\begin{remark} (i). Under the setting of Corollary~\ref{corr: case2.asymp},  the following asymptotic convergence result holds,
\[
\lim_{K \to +\infty} \mathbb{E} [ \text{dist}_{ \| \cdot \|_\infty} (0, \partial(f+\delta_X)(\x^{Y+1})) ]  = 0.
\]
Therefore, 
\[
\liminf_{k \to +\infty} \mathbb{E} [ \text{dist}_{ \| \cdot \|_\infty} (0, \partial(f+\delta_X)(\x^{k+1})) ]  = 0.
\]
(ii).  Under the setting of Corollary~\ref{corr: case2.asymp} and $\epsilon  \in (0,1)$,  the algorithm will output a point such that 
$$ \mathbb{E} [ \text{dist}_{ \| \cdot \|_\infty} (0, \partial(f+\delta_X)(\x^{Y+1})) ] \le C_2 \epsilon $$ where $C_2$ is an absolute constant if
$$
m_1 = e^2 (\log d)\sigma_\infty^2,  \; \frac{\log K + 1}{\sqrt{K+1} - 1} \le \epsilon,  \;  \frac{\Delta_f L}{ \sqrt{K+1} - 1} \le \epsilon.  
$$
It can be verified that it suffices to let $K \textcolor{black}{\triangleq} \tilde{C}( \epsilon^{-2} \max\{ ( \log(\epsilon^{-1}) )^2, \Delta_f^2 L^2, 1 \} ) $ where $\tilde{C}$ is a sufficiently large absolute constant. Then the sample complexity guarantee under this setting is
\begin{align*}
\sum_{k=1}^K m_k = \sO( m_1 K^2 ) = \sO \left( ( \log d ) \sigma_\infty^2 \epsilon^{-4}  \max\{ ( \log(\epsilon^{-1}) )^4, \Delta_f^4 L^4 ,1 \} \right).
\end{align*}
\qed
\end{remark}
}
Apply variance reduction to case 2, and the following statements hold.
\begin{theorem}\label{thm: case2svrg}
For Algorithm~\ref{alg: DI-SGD} with \eqref{ineq: proj}, consider {\revv case 2 \eqref{case2phi}} and estimating $G^k$ using variance reduction {\revv \textcolor{black}{as in \eqref{def: Gk-SPIDER} and \eqref{mk_spider}.} Suppose that Assumption~\ref{assumption first-order oracle},~\ref{ass: subG},~\ref{ass: Lip},~\ref{ass: setw3hj} hold.} Let $\eta_k \equiv \eta > 0$,  $\psi_k \equiv \psi > 0$,  $p_k \equiv 1/K$. Then the following holds when $t_k \equiv t > 0$.
\begin{align}
\notag
    & \bE[ {\rm dist}_{\| \cdot \|_\infty}(0, \partial(f+\delta_X)(\x^{Y+1})) ] \\ 
  \notag
& \le \left( L + \frac{1}{\eta} + \sqrt{\frac{4\Omega^2 L^2 q_0}{m}} \right) \psi + \sqrt{ \frac{\Omega e^2 (\log d)\sigma_\infty^2}{m_1} }  \\ 
    \label{errbd10}
    & + \frac{1}{\psi} \left( \frac{\Delta_f}{K} + \frac{ t \Omega e^2 (\log d)\sigmainf^2}{2 m_1} + \left(\frac{2 t\Omega^2 L^2 q_0}{m} + \frac{L}{2} + \frac{1}{2t} \right) \psi^2 \right).
\end{align}
\end{theorem}
\begin{proof}
    See Appendix~\ref{app: proof}. \qed
\end{proof}
\begin{corollary}\label{corr: vr.comp.case2}
Under the settings of Theorem~\ref{thm: case2svrg}, let
\begin{align*}
    & \eta \triangleq \frac{1}{L},\; K \triangleq \left\lceil \frac{\Delta_f L}{\epsilon^2} \right\rceil,\; m_1 \triangleq \left\lceil \frac{ \Omega e^2 (\log d) \sigma_\infty^2}{\epsilon^2} \right\rceil,\; m \triangleq \Omega^2 q_0, \; mq_0 = m_1,\; \\
    & t \triangleq \frac{1}{L},\; \psi \triangleq \frac{\epsilon}{L}.
\end{align*}
Then by \eqref{errbd10} we have
\begin{align*}
\bE[ {\rm dist}_{\| \cdot \|_\infty}(0, \partial(f+\delta_X)(\x^{Y+1})) ] {\rev \le} C_3 \epsilon,
\end{align*}
where $C_3$ is an absolute constant independent of any parameter of the algorithm or the problem data. The sample complexity is upper bounded by
\begin{align*}
  {\revi \sO\left( \frac{ \Delta_f L (\log d)^2 \sigma_\infty}{\epsilon^3} \right).}
\end{align*}
\end{corollary}
{\revv \begin{remark}
    The sample complexity results in Corollary~\ref{corr: mb.comp.case1},~\ref{corr: mb.comp.case2} and Corollary~\ref{corr: vr.comp.case1},~\ref{corr: vr.comp.case2} almost match the lower bounds in their corresponding settings up to $\mbox{poly}(\log d)$ factors. We discuss details of the lower bounds in Appendix~\ref{app: lb}.
\end{remark}}

\section{Calculation of the subproblem}\label{sec: subprob}
In this section, we discuss the formula to calculate the proximal projection subproblem \eqref{def: proj} for various feasible regions including the unconstrained case, box constraints, and $\ell_1$-ball + box constraints. General polyhedral $X$ is also discussed. Note that the objective function is non-smooth, so general solvers such as (projected) gradient descent may not enjoy fast linear convergence, not to mention feasible regions without easy projection formula. We mainly focus on the first choice of $\phi$ in this section. Calculation for the second case will be analogous.
\paragraph{Case 1: $\phi(\x) \textcolor{black}{\triangleq} \frac{\hat \rho}{2} \| \x \|_1^2$}
\paragraph{Unconstrained.} Suppose that $X = \bR^d$. Then for such a choice of $\phi$, \eqref{def: proj} is equivalent to the following:
\begin{align}\label{opt: l1squareprox}
   \min_{\z} \quad \frac{1}{2} \| \z - \vbf \|^2 + \frac{\hat \rho}{2} \| \z \|_1^2.
\end{align}
The next lemma demonstrates that it is indeed tractable to find the solution of \eqref{opt: l1squareprox}.
\begin{lemma}\label{lm: comp.l12prox}
Suppose that $\z^*$ denotes the optimal solution of \eqref{opt: l1squareprox}. Then, if $\vbf = 0$, $\z^* = 0$. Otherwise, suppose that $\vbf \neq 0$. Let $i_1, i_2, \hdots, i_d$ be the permutation of $\{1,\hdots, d \}$ such that
\begin{align*}
    |\vbf_{i_1}| \le |\vbf_{i_2}| \le \hdots \le |\vbf_{i_d}|.
\end{align*}
Denote a dummy index $i_0 \textcolor{black}{\triangleq} 0$ and dummy scalar $\vbf_{i_0} \textcolor{black}{\triangleq} 0$. Let $s_k \triangleq \sum_{t=0}^{k-1} | \vbf_{i_t} |$ $ + \frac{(d-k+1)\hrho + 1}{\hrho} | \vbf_{i_k } |$, $\forall k = 0,...,d$. (By convention, $\sum_{t=0}^{-1} | \vbf_{i_t} | \triangleq 0$). Then $s_k$ is non-decreasing and $\| \vbf \|_1 < s_d$. Suppose that $\bar k$ satisfies $s_{\bar k} \le \| \vbf \|_1 < s_{\bar k + 1}$. Then we have that
\begin{align}\label{proxform2}
    \z_{i_t}^* = \begin{cases}
        0 & \mbox{if } 1 \le t \le \bar k \\
        \vbf_{i_t} - \frac{\sgn(\vbf_{i_t}) \hrho}{\hrho(d-k)+1} \sum_{t = k+1}^d | \vbf_{i_t} | & \mbox{if } t > \bar k.
    \end{cases}
\end{align}
\end{lemma}
\begin{proof}
When $\vbf = 0$, it is trivial to see that $\z^* = 0$. Suppose that $\vbf \neq 0$. Note that we have 
\begin{align}
\notag
  0 & \in \z^* - \vbf + \frac{\hrho}{2} \partial \| \z \|_1^2 \mid_{\z = \z^*} \overset{\eqref{l1subdiff}}{=} \z^* - \vbf + \hrho  \| \z^* \|_1 \partial \| \z \|_1 \mid_{\z = \z^*} \\
  \notag
  \Longleftrightarrow \quad \z^* & = \mbox{arg}\min_{\z} \frac{1}{2} \| \z - \vbf \|^2 + \hrho \| \z^* \|_1 \| \z \|_1 \\
  \label{proxform1}
  \Longleftrightarrow \quad \z_i^* & = \begin{cases}
      0 & \mbox{if } | \vbf_i | \le \hrho \| \z^* \|_1 \\
      \vbf_i - \hrho \| \z^* \|_1 & \mbox{if } \vbf_i > \hrho \| \z^* \|_1 \\ 
      \vbf_i + \hrho \| \z^* \|_1 & \mbox{if } \vbf_i < -\hrho \| \z^* \|_1
  \end{cases}
\end{align}
It is easy to see that $\hrho \| \z^* \|_1 \le | \vbf_{i_d} |$, otherwise by \eqref{proxform1} $ \z^* = 0$, a contradiction. If $\hrho \| \z^* \|_1 = | \vbf_{i_d} |$, then by \eqref{proxform1}, $\| \z^* \|_1 = | \vbf_{i_d} | = 0 \implies \vbf = 0$, a contradiction to our assumption. Therefore, $\hrho \| \z^* \|_1 < | \vbf_{i_d} |$. Suppose that $ | \vbf_{i_k} | \le \hrho \| \z^* \|_1 < | \vbf_{i_{k+1}} |$ for some $k$ s.t. $0 \le k \le d-1$. Note that such $k$ is unique for a fixed $\vbf$. Then by \eqref{proxform1} we have
\begin{align*}
    \z_{i_t}^* = \begin{cases}
        0 & \mbox{if } 1 \le t \le k \\
        \vbf_{i_t} - \sgn(\vbf_{i_t}) \hrho \| \z^* \|_1 & \mbox{if } t > k
    \end{cases}
\end{align*}
This indicates that
\begin{align*}
    & \| \z^* \|_1 = \sum_{t = k+1}^d | \vbf_{i_t} | - \hrho (d - k) \| \z^* \|_1 \implies \| \z^* \|_1 = \frac{1}{\hrho(d-k) +1} \sum_{t = k+1}^d | \vbf_{i_t} |.
\end{align*}
Therefore,
\begin{align}
    \z_{i_t}^* = \begin{cases}
        0 & \mbox{if } 1 \le t \le k \\
        \vbf_{i_t} - \frac{\sgn(\vbf_{i_t}) \hrho}{\hrho(d-k)+1} \sum_{t = k+1}^d | \vbf_{i_t} | & \mbox{if } t > k
    \end{cases}
\end{align}
By the fact that $ | \vbf_{i_k} | \le \hrho \| \z^* \|_1 < | \vbf_{i_{k+1}} | $, this solution leads to $s_k \le \| \vbf \|_1 < s_{k+1}$. By the fact that $s_{t+1} - s_t = \frac{(d-t) + 1}{\hrho} (| \vbf_{i_{t+1}}  | - | \vbf_{i_t} |)  $ and $s_d = \| \vbf \|_1 + \frac{1}{\hrho}| \vbf_{i_d} | $, $s_t$ is non-decreasing and $\| \vbf \|_1 < s_d$. \qed
\end{proof}
\begin{remark}
    According to the formula given by Lemma~\ref{lm: comp.l12prox}, exact computation of $\z^*$ of \eqref{opt: l1squareprox} only requires $\sO(d)$ fundamental operations. Resolution of \eqref{opt: l1squareprox} is also discussed in \cite[Lemma 6.70]{beck2017first}, providing a different form of formula necessitating solving a one-dimensional root-finding problem.
\end{remark}
\paragraph{Box constraints.} {\revv Then we discuss the solution to the proximal projection problem \eqref{def: proj} when $X$ includes box constraints. In this case, we aim to solve a problem as follows:
\begin{align}\label{opt: l1squareprox-box}
   \min_{\z} \quad \frac{1}{2} \| \z - \vbf \|^2 + \frac{\hat \rho}{2} \| \z \|_1^2 \quad \mbox{subject to} \quad l \le \z \le u
\end{align}
where $\vbf \in\R^{n}$ and the box constraints $-\infty<l_i<u_i<+\infty$.

Introduce the shrinkage level \begin{align}\label{eq:def_of_tau}
    \tau\textcolor{black}{\triangleq}\hat\rho\|\bm z\|_{1}\;(>0)
\end{align} and dual variables $\lambda_i^{\mathrm L},\lambda_i^{\mathrm U}\ge0$ for the lower and upper bounds. Then the KKT conditions of such problem is for each $i$,
\begin{align}
    0 &= \z_i-\vbf_i+\tau\,\xi_i+\lambda_i^{\mathrm U}-\lambda_i^{\mathrm L}\nonumber\\
  \xi_i &\in\partial|\z_i|=
        \begin{cases}
           \{\sgn(\z_i)\}, & \z_i\neq0\\
           [-1,1],                       & \z_i=0
        \end{cases}\nonumber\\
  0 &\le\lambda_i^{L}\;\perp\;\z_i-l_i\ge0\nonumber\\
  0 &\le\lambda_i^{U}\;\perp\;u_i-\z_i\ge0\nonumber
\end{align}

\begin{itemize}
    \item If $\z_i\in(l_i,u_i)$: we find that $\xi_i=\sgn(\z_i)$ when $\z_i \neq 0$, so we could obtain $\sgn(\z_i)=\sgn(\vbf_i)$ and $\z_i=\vbf_i-\tau\sgn(\z_i)$. If $\z_i=0$, $\lvert \vbf_i\rvert<\tau$. Combine them and we get \begin{align*}
    \z_i
  =\sgn(\vbf_i)\bigl(\lvert \vbf_i\rvert-\tau\bigr)_+.
\end{align*}

  \item If $\z_i\in\{l_i, u_i\}$: then $\sgn(\vbf_i)(| \vbf_i |-\tau)_+$ would fall outside the interval $(l_i, u_i)$. Therefore a nonnegative $\lambda_i^L$ or $\lambda_i^U$ will compensate the imbalance in the first equation of KKT conditions caused by the clipping operation as follows. \begin{align*}
      \z_i = \operatorname{clip}_{[l_i, u_i]}\left(\sgn(\vbf_i)(\lvert \vbf_i\rvert - \tau)_+\right).
  \end{align*}
  
\end{itemize}
So in general, \begin{align}\label{eq:find_z}
    \z_i = \operatorname{clip}_{[l_i, u_i]}\left(\sgn(\vbf_i)(\lvert \vbf_i\rvert - \tau)_+\right).
\end{align}
Substitute it back to \eqref{eq:def_of_tau} and we get the equation \begin{align}
    \tau =\hat\rho\sum_{i=1}^{n}\left\vert
        \operatorname{clip}_{[l_i,u_i]}
           \left(\operatorname{sign}(\vbf_i)(\lvert \vbf_i\rvert-\tau)_+\right)\right\vert.
\end{align}
So the above problem is transformed to a root finding problem of such monotone function \begin{align}\label{def: res.search}
    \sR(\tau) \triangleq \tau - \hat\rho\sum_{i=1}^{n}\left\vert
        \operatorname{clip}_{[l_i,u_i]}
           \left(\operatorname{sign}(\vbf_i)(\lvert \vbf_i\rvert-\tau)_+\right)\right\vert.
\end{align}
Also note that $\sR(0)\le 0$ and $\sR(\hrho \sum_{i=1}^d\max\{|\vbf_i|,|l_i|,|u_i|\})\ge 0$. Therefore we could use bisection to search the unique solution $\tau^*$ in the interval \\
$[0, \sum_{i=1}^d\max\{|\vbf_i|,|l_i|,|u_i|\}]$, and with that we could obtain the solution of the subproblem \eqref{opt: l1squareprox-box} from
\eqref{eq:find_z}.

\paragraph{$\ell_1$ ball + box constraints.} We make the feasible more complicated by incorporating both an $\ell_1$ ball and box constraints. It is not obvious how to efficiently calculate the projection onto this feasible region. In this case, \eqref{def: proj} is equivalent to the following subproblem:

\begin{align}\label{opt: l1squareprox-boxl1}
\begin{array}{rl}
    \min\limits_{\x} & \quad \frac{1}{2} \| \x - \vbf \|^2 + \frac{\hrho}{2} \| \x \|_1^2 \\  \mbox{subject to} & \quad \| \x - \w \|_1 \le \alpha, \; l \le \x \le u. 
\end{array}
\end{align}
We first write the KKT conditions for this problem:
\begin{align}
\label{eq1}
    & 0 = \x - \vbf + \hrho \| \x \|_1 \xi + \mu \xi^w - \lambda^L + \lambda^U, \\
    \notag
    \xi_i &\in\partial|\x_i|=
        \begin{cases}
           \{\sgn(\x_i)\}, & \x_i\neq0\\
           [-1,1],                       & \x_i=0
        \end{cases}\nonumber\\
    \notag
    \xi_i^w &\in\partial|\x_i - \w_i|=
        \begin{cases}
           \{\sgn(\x_i - \w_i)\}, & \x_i\neq \w_i \\
           [-1,1],                       & \x_i= \w_i
        \end{cases}\nonumber\\
    \notag
    & 0 \le \mu \perp \| \x - \w \|_1 - \alpha \le 0, \\
    \label{eq1.3}
    & 0 \le \lambda^L \perp l - \x \le 0, \quad 0 \le \lambda^U \perp \x - u \le 0.
\end{align}
Our goal is to solve the above equations. 
First, we suppose that $\mu = 0$ and let $\tau \textcolor{black}{\triangleq} \hrho \| \x \|_1$ and treat it as a free variable. Then \eqref{eq1} is equivalent to 
\begin{align}
\label{eq2}
    0 = \x - \vbf + \tau \xi - \lambda^L + \lambda^U.
\end{align}
Solution to \eqref{eq2} is
\begin{align*}
    \bar \x_i(\tau) = \mbox{clip}_{[l_i,u_i]}( \sgn(\vbf_i)(\vbf_i - \tau)_+ ).
\end{align*}
Then we try to solve 
\begin{align}\label{def: res.search1}
\sR_1(\tau) \triangleq \tau - \hrho \sum_{i} | \mbox{clip}_{[l_i,u_i]}( \sgn(\vbf_i)(\vbf_i - \tau)_+ ) | = 0
\end{align}
by searching $\tau$ via bisection. Note that $\sR_1(\tau)$ is continuous and monotonically nondecreasing with $\sR_1(0) \le 0$ and $\sR_1( \hrho \sum_i \max\{ |l_i|, |u_i| \} ) \ge 0$. Denote the root of $\sR_1$ as $\tau^*$. Then check if
\begin{align*}
    \| \bar \x(\tau^*) - \w \|_1 \le \alpha.
\end{align*}
If yes, stop the algorithm and output $\bar \x(\tau^*)$. Otherwise, we need to determine the correct $\mu$. Now suppose that we are given a fixed $\mu > 0$ and $\w_i$. By observation of \eqref{eq1} and \eqref{eq1.3}, and letting $\tau \textcolor{black}{\triangleq} \hat \rho \| \x \|_1$, we have the following for solution $\bar \x$ that satisfies \eqref{eq1} and \eqref{eq1.3} (suppose that $\w_i \neq 0$, otherwise we have similar formula as in the case $\mu = 0$).
\begin{align*}
    \bar \x_i = \begin{cases}
    \mbox{clip}_{[l_i,u_i]} (\tau + \vbf_i + \mu) & \mbox{if } \tau < -\vbf_i - \mu + \min(0,\w_i) \\
    \mbox{clip}_{[l_i,u_i]} (-\tau + \vbf_i - \mu) & \mbox{if } \tau < \vbf_i - \mu - \max(0,\w_i) \\
    \mbox{clip}_{[l_i,u_i]} (-\tau \sgn(\w_i) + \vbf_i + \mu \sgn(\w_i)) & \begin{array}{l}
         \mbox{if } \tau \le \vbf_i\sgn(\w_i) + \mu \\
         \mbox{\& } \tau \ge \vbf_i\sgn(\w_i) + \mu - |\w_i|
    \end{array}\\
    \mbox{clip}_{[l_i,u_i]}(0) & \mbox{if } \tau \ge | \vbf_i + \mu \sgn(\w_i) | \\
    \mbox{clip}_{[l_i,u_i]} (\w_i) & \begin{array}{l}
         \mbox{if } \tau \le \vbf_i\sgn(\w_i) + \mu - |\w_i| \\
         \mbox{\& } \tau \ge \vbf_i\sgn(\w_i) - \mu - |\w_i|
    \end{array}
    \end{cases}
\end{align*}
Alternatively, we have the following equivalent formula,
\begin{align*}
    \bar \x_i =
        \mbox{clip}_{[l_i,u_i]}( \tau \sgn(\w_i) + \vbf_i  + \mu \sgn(\w_i) ) \quad \mbox{when } - \vbf_i \sgn(\w_i) - \mu \ge 0.
\end{align*}
When $- \vbf_i \sgn(\w_i) - \mu < 0$, we have that 
\begin{align*}
    & \bar \x_i = \\
    & \begin{cases}
        \mbox{clip}_{[l_i,u_i]}( 0 ) & \mbox{when }
        \tau \ge \vbf_i \sgn(\w_i) + \mu \\
        \mbox{clip}_{[l_i,u_i]}( -\tau \sgn(\w_i) + \vbf_i  + \mu \sgn(\w_i) ) &  \mbox{when } \tau \ge \vbf_i \sgn(\w_i) + \mu - |\w_i| \\
        & \mbox{and } \tau \le \vbf_i \sgn(\w_i) + \mu \\
        \mbox{clip}_{[l_i,u_i]}( \w_i ) & \mbox{when } \tau \ge \vbf_i \sgn(\w_i) - \mu - |\w_i| \\
        & \mbox{and } \tau \le \vbf_i \sgn(\w_i) + \mu - |\w_i| \\
        \mbox{clip}_{[l_i,u_i]}(-\tau\sgn(\w_i) + \vbf_i - \mu\sgn(\w_i) ) & \mbox{when } \tau < \vbf_i \sgn(\w_i) - \mu - |\w_i|
    \end{cases}
\end{align*}
Note that $|\bar \x_i(\tau)|$ is monotonically non-increasing. Then we try to solve
\begin{align}\label{def: res.search2}
    \sR_2(\tau) \triangleq \tau - \hrho \sum_i | \bar \x_i(\tau) | = 0
\end{align}
via searching $\tau$. $\sR_2(\tau)$ is continuous and monotonically nondecreasing with $\sR_2(0) \le 0$ and $\sR_2(\hrho\sum_i M_i) \ge 0$, where $M_i \triangleq \max\{ |l_i|, |u_i| \}$. Therefore, the solution to \eqref{def: res.search2} exists.
We denote $\bar \x(\mu, \tau^*(\mu))$ as the solution when $\mu$ is given. Search $\mu > 0$ such that 
\begin{align}\label{res.search3}
    \sR_3(\mu) \triangleq \| \bar \x(\mu, \tau^*(\mu)) - \w \|_1 - \alpha = 0.
\end{align}
Note that $\sR_3(0) > 0$ (otherwise $\mu = 0$ and we conclude the discussion). $\sR_3(\mu)$ is continuous and montonically nonincreasing (In fact, $\mu$ can be viewed as the coefficient of the regularizer $ \| \x - \w \|_1 $). \\
Recall that 
\begin{align*}
    & 0 = \x - \vbf + \hrho \| \x \|_1 \xi + \mu \xi^w - \lambda^L + \lambda^U, \\
    & 0 \le \lambda^L \perp l - \x \le 0, \quad 0 \le \lambda^U \perp \x - u \le 0.
\end{align*}
Let $\hat \x \textcolor{black}{\triangleq} \mbox{clip}_{[l,u]}(\w)$ and $\hat \tau \textcolor{black}{\triangleq} \hrho \| \hat \x \|_1$. Then if we let $\hat \mu \textcolor{black}{\triangleq} \| \hat \x - \vbf \|_\infty + \hat \tau$, $\hat \mu$ and $\hat x$ satisfy the above conditions for certain choice of $\lambda^L, \lambda^U, \xi, \xi^w$. In this case, $\| \hat x - \w \|_1 - \alpha \le 0$ otherwise \eqref{opt: l1squareprox-boxl1} is not feasible. Therefore, 
$
    \sR_3(\hat \mu) \le 0.
$
Solution to \eqref{res.search3} exists.

\paragraph{Polyhedron.} Last but not least, we consider resolution of the following general subproblem:
\begin{align}\label{opt: l1squareprox-affine}
\begin{array}{rl}
    \min\limits_{\x} & \quad \frac{1}{2} \| \x - \vbf \|^2 + \frac{\hrho}{2} \| \x \|_1^2 \\  \mbox{subject to} & \quad A \x \ge b
\end{array}
\end{align}
By splitting $\x$ into positive and negative parts $\x = \x^+ - \x^-$ and letting $\z \textcolor{black}{\triangleq} \begin{pmatrix} \x^+ \\ \x^- \end{pmatrix}$, the above subproblem can be rewritten as
\begin{align*}
    \begin{array}{rl}
    \min\limits_{\z} & \quad \frac{1}{2} \z^T Q \z + \mathbf{p}^T \z \\  \mbox{subject to} & \quad \tilde{A} \z \ge b, \; \z \ge 0.
\end{array}
\end{align*}
where $Q \textcolor{black}{\triangleq} \begin{pmatrix} I & -I \\ -I & I \end{pmatrix} + \hrho \begin{pmatrix} 1 & \hdots & 1 \\ \vdots & \ddots & \vdots \\ 1 & \hdots & 1
\end{pmatrix}$, $\mathbf{p} \textcolor{black}{\triangleq} \begin{pmatrix} -\vbf \\ \vbf \end{pmatrix}$, $\tilde{A} \textcolor{black}{\triangleq} \begin{pmatrix} A, & -A \end{pmatrix}$. i.e., \eqref{opt: l1squareprox-affine} is equivalent to a convex quadratic program (QP) and we can apply efficient QP solvers for resolution.
}
\paragraph{Case 2: $\phi(\x) \textcolor{black}{\triangleq} \delta_{\{ \z \mid \| \z \|_1 \le \psi\}}(\x)$}
\paragraph{Unconstrained.}
In this case, we need to solve the following problem in the proximal projection step.
\begin{align}\label{opt: l1consprox}
\begin{aligned}
\min_{\z} \quad &\dfrac{1}{2} \lVert \z-\vbf \rVert_2^2 \\
\textit{s.t.} \quad &\lVert \z \rVert_1 \leq \psi
\end{aligned}
\end{align}
The next lemma reveals the closed-form solution to \eqref{opt: l1consprox}. Since this result utilizes similar technique as in Lemma~\ref{lm: comp.l12prox}, we relegate its proof to Appendix~\ref{app: proof}.
\begin{lemma}\label{lm: l1consprox}
Suppose that $\z^*$ denotes the optimal solution of \eqref{opt: l1consprox}. If $\| \vbf \|_1 \le \psi$, then $ \z^* = \vbf $. Otherwise, suppose that $i_1,\hdots,i_d$ is a permutation of $\{1,2,\textcolor{black}{\ldots},d \}$ such that $ | \vbf_{i_1}| \ge | \vbf_{i_2}| \ge \textcolor{black}{\ldots} \ge | \vbf_{i_d}|$. Denote $s_0 \triangleq 0$, $s_d \triangleq \| \vbf \|_1$ and $s_m \triangleq \sum_{i=1}^m | \vbf_{i_m} | - m | \vbf_{i_{m+1}} |$, $m = 1,2,\hdots,d-1$. Then $s_m$ is a non-decreasing sequence. Suppose that $s_{\bar m-1} < \psi \le s_{\bar m}$ for some $\bar m \in [d]$. Then we have
\begin{align*}
    & \z_{i_t}^* = \begin{cases}
        0 & \mbox{if }  \bar m < t \le d \\
        \vbf_{(i_t)} - \frac{\sgn\left(\vbf_{(i_t)}\right)}{\bar m}\left(\sum_{k=1}^{\bar m} \left\lvert \vbf_{(i_k)} \right\rvert - \psi \right) & \mbox{if }  1 \le t \le \bar m
    \end{cases}
\end{align*}
\end{lemma}
\begin{remark}
    According to the formula in Lemma~\ref{lm: l1consprox}, exact computation of $z^*$ only requires $O(d)$ fundamental operations.
\end{remark}
{\revv
\begin{remark}
    When the feasible region is box constraints, $\ell_1$ ball + box constraints, and general polyhedra, the solution of the subproblem in case 2 is in similar fashion to that in case 1 so we omit the discussion.
\end{remark}}

\section{Numerical experiment}\label{section:num}
We conduct several experiments in this section to illustrate our findings. In Section~\ref{subsec: subprob}, we compare our formula for solving the subproblems with a general solver - ADMM and demonstrate the efficiency of the former. In Section~\ref{sec:exp_1}, we run our algorithms to solve \eqref{SP problem} and illustrate their dimension insensitive property. The numerical experiments in this section are carried on a \textit{32 vCPU Intel(R) Xeon(R) Platinum 8352V CPU @ 2.10GHz}. 
\subsection{Efficiency of the specialized solvers for subproblems}\label{subsec: subprob}
In this section, we compare the solvers proposed in Section~\ref{sec: subprob} with a general constrained subproblem solver ADMM \cite{boyd2011distributed}. We test the running time of different solvers on the subproblem with only the box constraints \eqref{opt: l1squareprox-box}, and with $\ell_1$ ball + box constraints \eqref{opt: l1squareprox-boxl1}. Details of the solution methods are as follows. 
\begin{enumerate}
    \item \textbf{One/two-dimensional search}. As we elaborated in Section~\ref{sec: DISFOM}, the subproblem with box constraint \eqref{opt: l1squareprox-box} is equivalent to finding the root of a one-dimensional monotone function. We use bisection to search the root 
    The subproblem with $\ell_1$-ball + box constraint \eqref{opt: l1squareprox-boxl1} is equivalent to root finding of two nested one-dimensional function. We use bisection to search the roots.
    \item \textbf{ADMM}. In order to use ADMM \cite{yang2016linear}, we reformulate our problem as follows: 
    \begin{align*}
\begin{array}{rl}
    \min\limits_{\x \in \bR^{d_1},\y \in \bR^{d_2}} & \quad g(\x) + h(\y) \\  \mbox{subject to} & \quad A \x + B \y = b. 
\end{array}
\end{align*}
For the box constrained subproblem \eqref{opt: l1squareprox-box}, $d_1 = d_2$, $g(\x) = \frac{1}{2} \| \x - \vbf \|^2 + \frac{\hrho}{2} \| \x \|_1^2 $, $h(\y) = \delta_{l \le \y \le u} (\y)$, $A = I_{d_1}$, $B = - I_{d_1}$, $b = 0$. For the $l_1$-ball + box constrained subproblem \eqref{opt: l1squareprox-boxl1}, $2d_1 =  d_2$, $g(\x) = \frac{1}{2} \| \x - \vbf \|^2 + \frac{\hrho}{2} \| \x \|_1^2$, $ \y = [\y_{(1)}^\top,\y_{(2)}^\top]^\top \in \bR^{d_1} $, $\y_{(1)}, \y_{(2)} \in \bR^{d_1}$, $ h(\y) = \delta_{\| \y \|_1 \le \alpha} (\y_{(1)}) + \delta_{l \le \y \le u} (\y_{(2)})  $, $A = [I_{d_1},I_{d_1}]^\top$, $B = -I_{d_2}$, $b = [\w;0]$, $\w \in \bR^{d_1}$. Through such decomposition, each iteration of ADMM can be solved in closed form and the convergence rate is linear according to \cite{yang2016linear}. ADMM is a natural and efficient choice among general solvers to tackle the subproblem because it exploits the problem structure and enjoys fast convergence.
\end{enumerate}
     
\begin{figure}[ht]
    \centering
    \includegraphics[width = 0.49\textwidth]{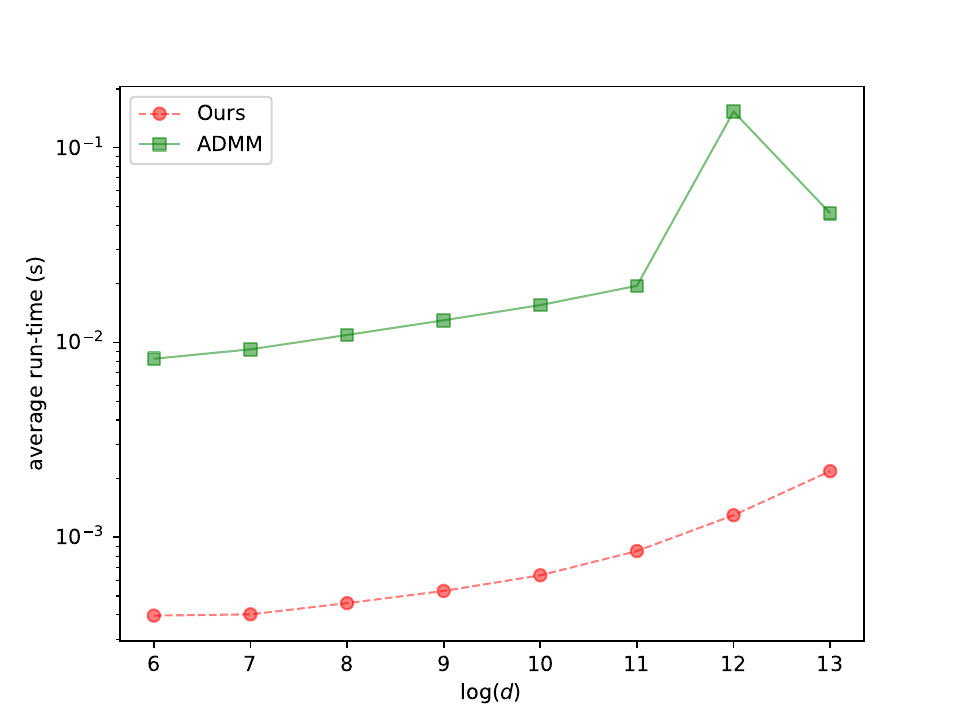}
    \includegraphics[width = 0.49\textwidth]{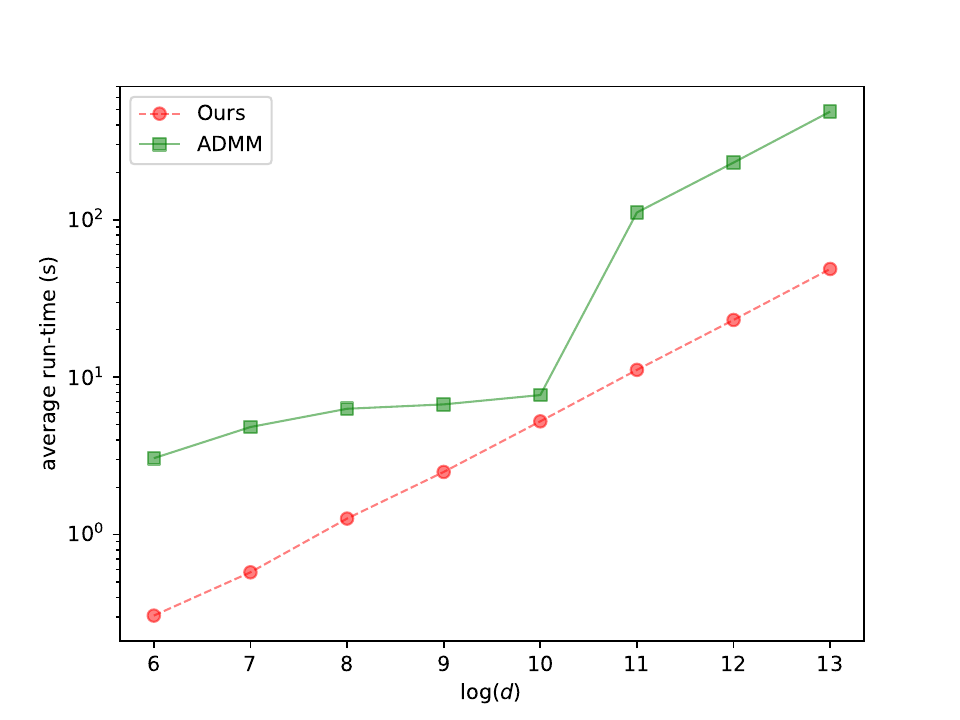}
    \caption{The average running time (over 10 trials) of our solvers in Section \ref{sec: subprob} and ADMM in Section \ref{subsec: subprob} in solving the subproblem with \textbf{only the box constraints} (left) and \textbf{$\ell_1$ ball + box constraints} (right).}
    \label{fig:time_test}
\end{figure}

Next we specify the details of other experiment settings as below. \begin{itemize}
    \item \textbf{Box constraints}: \quad For the subproblem \eqref{opt: l1squareprox-box}, we sample $\vbf, \wbf\in\mathcal{N}(0, I_d)$ and set $l_i \textcolor{black}{\triangleq} -20, u_i\textcolor{black}{\triangleq}20, \forall i$. We carry on our experiments with dimension $d\in\{2^6, 2^7, 2^8,\dots, 2^{13}\}$. $\hat \rho \textcolor{black}{\triangleq} 1$. For dimension $d$ we set the scaled ADMM penalty parameter $\beta \textcolor{black}{\triangleq} 0.1+0.3\log d$. We terminate our solver when $upper\ bound \ - \ lower\ bound < 10^{-10}$ for the search process of $\tau$ and $| \sR(\tau) |< 10^{-10}$ ($\sR$ defined in \eqref{def: res.search}). Solution of this solver is denoted as $\x_{sol}$. Then we terminate the ADMM solver when it reaches $100$ times of the running time of our solver, or the function value of $\y^{k+1}$ is smaller than the function value of our solver. Fig \ref{fig:time_test} (left) reports the results.

    \item \textbf{$\ell_1$ ball + box constraints} \quad For the subproblem \eqref{opt: l1squareprox-boxl1}, {we sample $\vbf, \wbf$ uniformly on $[-50,50]^d$ and $[-20,20]^d$ respectively. Set $\alpha \textcolor{black}{\triangleq} 10$, $l_i \textcolor{black}{\triangleq} -20, u_i\textcolor{black}{\triangleq}20, \forall i$. $\hat \rho \textcolor{black}{\triangleq} 1$.} We carry on our experiments with dimension $d\in\{2^6,2^7,$ $2^8,\dots, 2^{13}\}$. For each dimension we set the scaled ADMM penalty parameter $\beta \textcolor{black}{\triangleq} 100+300\log d$. We terminate our solver when $upper\ bound \ - \ lower\ bound < 10^{-6}$ for both line search process of $\tau$ and $\mu$ and $\|\bar{\xbf}(\mu, \tau^*(\mu)) $ $- \wbf\| - \alpha< 10^{-12}$. Solution of this solver is denoted as $\x_{sol}$. Then we terminate the ADMM solver when it reaches $10$ times of the running time of our solver, or when it satisfies the following conditions: \begin{align*}
    \|\ybf_{(2)}^{k+1}-\wbf\|_1 - \alpha & < 10^{-12}\\
    f(\ybf_{(2)}^{k+1}) & < f(\xbf_{sol})\end{align*}
    See Fig \ref{fig:time_test} (right) for the results. 
\end{itemize}

For the problem with box constraints (left), our solver is consistently faster than ADMM by at least 20 times. For the “$\ell_{1}$ + box" constrained problem (right), our solver also performs better. For 
half of all the cases, ADMM reaches the time limit, being $10$ times slower than our solver. In the high dimensional cases, the feasibility tolerance is violated by ADMM, which causes more computation time.
Note that ADMM is among the best solvers for the subproblem (easy implementation, linear convergent). This demonstrate the efficiency of our proposed algorithms to solve the subproblems.

\subsection{Dimension insensitive property illustration}\label{sec:exp_1}
In this section, we present some numerical experiments to illustrate the dimension insensitive property of our proposed methods.\footnotemark \footnotetext{Please refer to \href{https://github.com/bjw010529/DISFOM.git}{https://github.com/bjw010529/DISFOM.git} for details of the experiments.} 

\textbf{Problem setting} \quad We consider a nonconvex stochastic quadratic programming problem as follows:
\begin{align}
    \min_{\x \in X} \quad f(\x) \triangleq \frac{1}{2} \mathbb{E}\left[({\bm \alpha}^T \x - b)^2\right] + \lambda\sum_{i=1}^d \frac{\x_i^2}{1 + \x_i^2},
\end{align}
where the feasible region $X \textcolor{black}{\triangleq} [-R, R]^d$, $(\bm{\alpha}, b)\in \R^d\times\R $ is a random pair that satisfies linear relationship $b = \bm{\alpha}^T \x_{\rm true} + w $. We generate the scaled covariance matrix for $\bm{\alpha}$ by the following procedure. First we generate a $d\times d$ identity matrix, and then replace its top-left $100\times100$ principal sub-matrix by $\Sigma_{\rm sub} \textcolor{black}{\triangleq} \bm{QDQ}^T \in \R^{100\times100}$, where $\bm{Q}$ consists of the orthonormal basis of a $100\times100$  matrix whose entries are i.i.d. uniformly distributed on $(0,1)$ and $\bm{D}$ is a diagonal matrix with each diagonal entry i.i.d. uniformly distributed on $(1,2)$. Let $\bm \alpha \textcolor{black}{\triangleq} \Sigma^{\frac{1}{2}} s$, where $s_i$, $i = 1,...,d$ are i.i.d and obey truncated standard normal distribution over $[-u,u]$. $w$ also obeys a truncated standard normal distribution over $[-u,u]$ and is independent of $\bm \alpha$. Therefore, the variance $\sigma^2$ of $s_i$ and $w$ has the following formula
\begin{align}
    \sigma^2 = 1 - \frac{\frac{2u}{\sqrt{2\pi}} \exp\left(-\frac{u^2}{2} \right) }{\Phi(u) - \Phi(-u)},
\end{align}
where $\Phi$ is the CDF of standard normal distribution. By construction, function $f$ has the following closed form:
\begin{align*}
    f(\x) = \frac{\sigma^2}{2}(\x - \x_{\rm true})^T \Sigma (\x - \x_{\rm true}) + \lambda\sum_{i=1}^d \dfrac{\x_i^2}{1 + \x_i^2} + \frac{\sigma^2}{2}.
\end{align*}
The Lipschitz constant $L$ in 2-norm (thus in $\ell_\infty$/$\ell_1$-norm) of $\nabla f$ is $\lambda_{\max} + 2\lambda$, where $\lambda_{\max}$ is the largest eigenvalue of $\sigma^2 \Sigma$, ranging from $\sigma^2$ to $2\sigma^2$ by construction. Therefore, $L$ has a fixed range independent of $d$. $f$ is nonconvex when the minimal eigenvalue of $\nabla^2 f$, $\lambda_{\min} - \lambda/2$ is negative, where $\lambda_{\min}$ denotes the minimal eigenvalue of $\sigma^2 \Sigma$ and equals to $\sigma^2$. The problem setting also satisfies other assumptions in Section~\ref{sec: Prelim}. The sampled gradient has the formula  $\nabla_{\x} F(\x,\alpha,w) = \alpha(\alpha^T(\x-\x_{\rm true})-w)) + \lambda \sum_{i=1}^d \frac{2\x_i}{(1 + \x_i^2)^2}$. Its Lipschitz constant is $L(\alpha) = \lambda_{\max}(\alpha \alpha^T) + 2 \lambda \le \| \alpha \|^2 + 2 \lambda \le M < +\infty$ since $\alpha$ has bounded support. This uniform bound on the Lipschitz constant indicates Assumption~\ref{ass: Lip}. The gradient noise is $(\alpha \alpha^T - \sigma^2 \Sigma )(\x - \x_{\rm ture}) - \alpha w$, which is unbiased and has bounded support given that $\alpha, \w$ have bounded support and $\x$ is in a bounded region. Therefore it conforms sub-Gaussian distribution componentwisely. 

The experiments present the quality of solutions, measured by the averaged gap $f - f^*$ and residual out of 3 replications (for each algorithm and each measure, we divide the data by the quantity at $d = 2^7$). The value of $f^*$ is given by using the projected gradient method with backtracking (see Algorithm \ref{alg: GD backtracking} in Appendix~\ref{app: alg}), to solve the correlated closed-form problem.
We carry on our experiments with dimension $d \in \{2^7, 2^8, \dots, 2^{14}\}$, {\rev $\x_{\rm true} \textcolor{black}{\triangleq} (1,...,1,0,...,0)^T$ (first $100$ elements are 1)} and start with $\xbf_1 \textcolor{black}{\triangleq} \bm{0}$. $R \textcolor{black}{\triangleq} 3$, $u \textcolor{black}{\triangleq} 3$, $\lambda \textcolor{black}{\triangleq} 2.5$. 

\textbf{Gradient generation} \quad According to the two estimating methods of gradient mentioned above, mini-batch and variance reduction, we compare three algorithms with both methods of gradient estimation. In the mini-batch setting, the batch size is $m_k \equiv m \textcolor{black}{\triangleq} 1000$ and algorithms stop with $K \textcolor{black}{\triangleq} 300$. In the variance reduction setting, we set $m \textcolor{black}{\triangleq} 1000$ and $q_0 \textcolor{black}{\triangleq} 9$. The batch size is $m_k \textcolor{black}{\triangleq} m$ when \textcolor{black}{$k \equiv 1 \pmod{q_0}$} and $m_k \textcolor{black}{\triangleq} m^{\frac{2}{3}} = 100$ otherwise, and the algorithms stop with $K \textcolor{black}{\triangleq} 1350$. 

\textbf{Methods} Here we present the detailed implementation of DISFOMs, proximal stochastic gradient descent (SGD), proximal version of stochastic path integrated differential estimator (SPIDER) 
and stochastic mirror descent (SMDs).\begin{enumerate}

\item \textbf{DISFOM$\_$minibatch} \quad Algorithm \ref{alg: DI-SGD} with $\phi(\zbf) \textcolor{black}{\triangleq} \dfrac{\hat{\rho}}{2} \lVert \zbf \rVert_1^2 $ (Case 1) and minibatch sampling \eqref{def: Gk-SGD}. We choose $\hat{\rho} \textcolor{black}{\triangleq} 2$ in minibatch.
The stepsize is $\eta_k \equiv \eta \textcolor{black}{\triangleq} \dfrac{1}{L}$. 

\item \textbf{DISFOM$\_$vr} \quad Algorithm \ref{alg: DI-SGD} with $\phi(\zbf)  \textcolor{black}{\triangleq} \dfrac{\hat{\rho}}{2} \lVert \zbf \rVert_1^2 $ (Case 1) and variance reduction \eqref{def: Gk-SPIDER}.  We choose $\hat{\rho} \textcolor{black}{\triangleq} 128$. 
The stepsize is $\eta_k \equiv \eta \textcolor{black}{\triangleq} \dfrac{1}{L}$. 

\item (Proximal) \textbf{SGD} \quad Algorithm \ref{alg: DI-SGD} with $P_X^k \equiv P_X$ (Euclidean projection on $X$). $G^k$ is estimated using minibatch \eqref{def: Gk-SGD}. The step size is $\eta_k \equiv \eta \textcolor{black}{\triangleq} \dfrac{1}{L}$.

\item (Proximal) \textbf{SPIDER} \quad Algorithm \ref{alg: DI-SGD} with $P_X^k \equiv P_X$ (Euclidean projection on $X$). $G^k$ is estimated using variance reduction \eqref{def: Gk-SPIDER}. The step size is $\eta_k \equiv \eta \textcolor{black}{\triangleq} \dfrac{1}{10L}$.\footnotemark \footnotetext{Here we need to decrease the stepsize of SPIDER to stabilize its performance.}

\item \textbf{SMD$\_$minibatch} \quad Algorithm \ref{alg: DI-SGD}. 
Solve the proximal projection problem in Step 2.2 as follows:
 \begin{align*}
     \xbf^{k+1} = \mbox{arg}\min_{\zbf\in X} \left\{\langle G^k, \zbf\rangle + \dfrac{1}{\alpha_k} D_{\omega}(\zbf, \xbf^k) \right\}
 \end{align*}
 where $D_{\omega}(x,y) \textcolor{black}{\triangleq} \omega(x) - \omega(y) - \langle\nabla \omega(y), x - y \rangle $. Let the distance generating function $\omega(x) \textcolor{black}{\triangleq} \frac{C}{2}\lVert x\rVert_p^2$ where $p \textcolor{black}{\triangleq} 1 + \dfrac{1}{\ln d}$ and $C \textcolor{black}{\triangleq} e^2\ln d$ (so that $\omega(x)$ is 1-strongly convex \textcolor{black}{with respect to} $\| \cdot \|_1$), see \cite{beck2003mirror} and the reference therein. The step size is $\alpha_k \equiv \alpha \textcolor{black}{\triangleq} \dfrac{c}{\sqrt{K}}$, where $c \textcolor{black}{\triangleq} \sqrt{\dfrac{f(x^1)}{\rho  L^2}}$, $\rho \textcolor{black}{\triangleq} \dfrac{\lambda}{2} - \lambda_{\min}$, given by \cite{zhang2018convergence}. $G^k$ is estimated via minibatch.

\item \textbf{SMD$\_$vr} \quad Same settings as in SMD$\_$minibatch, except that $G^k$ is estimated via variance reduction.

\end{enumerate}

\textbf{Results and interpretation}
\begin{figure}[ht]
    \centering
    \includegraphics[width = 0.45\textwidth]{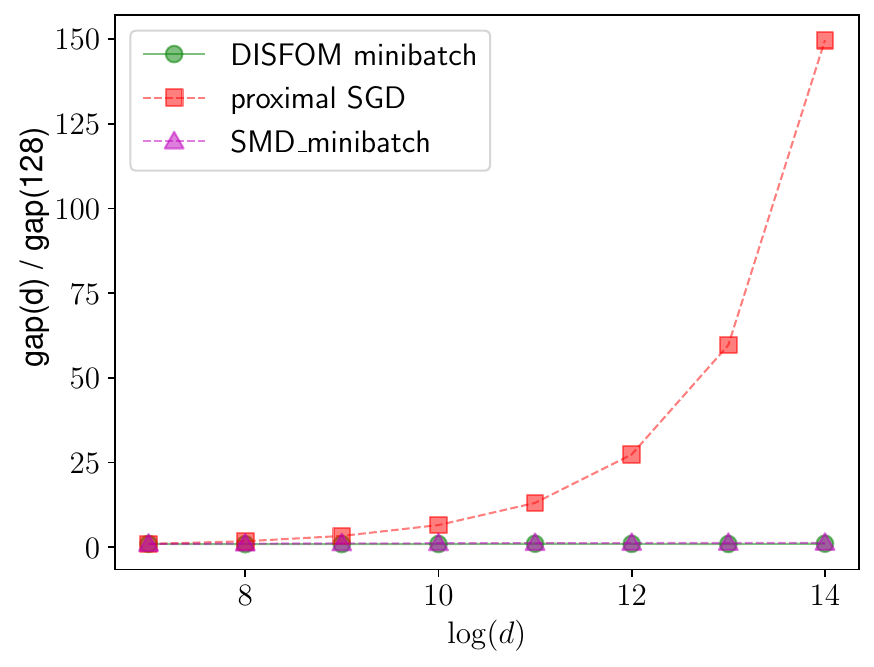}
    \includegraphics[width = 0.45\textwidth]
    {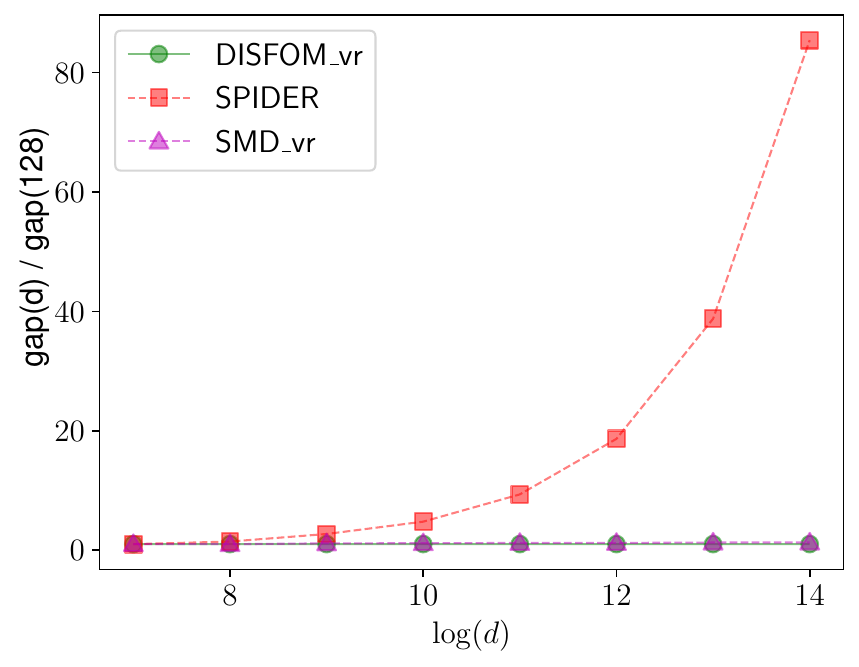}
    \caption{Averaged {relative} gap vs $\log d$ for minibatch (left) and variance reduction (right)
    }\label{fig:minibatch}
\end{figure}
\begin{figure}[ht]
    \centering
    \includegraphics[width = 0.45\textwidth]{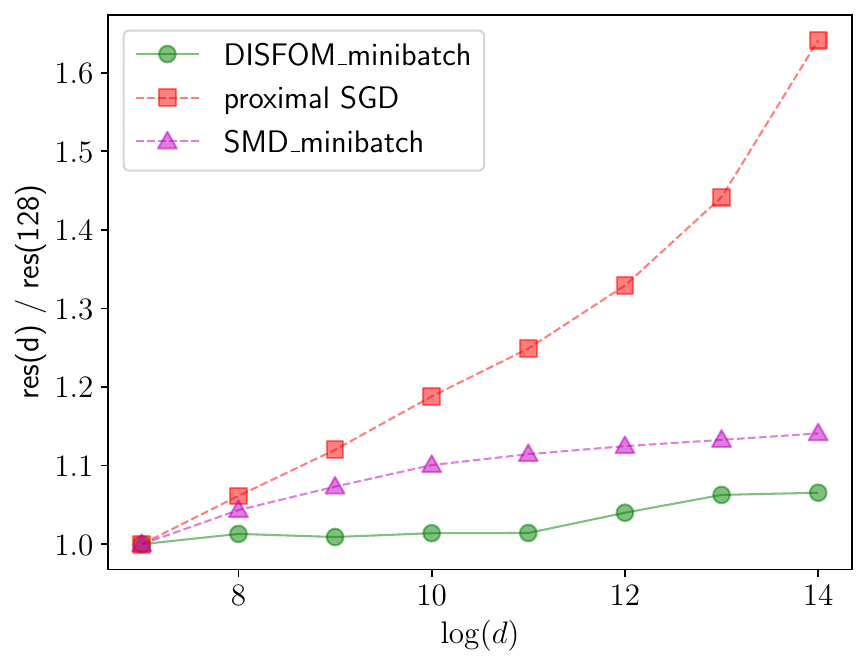}
    \includegraphics[width = 0.45\textwidth]{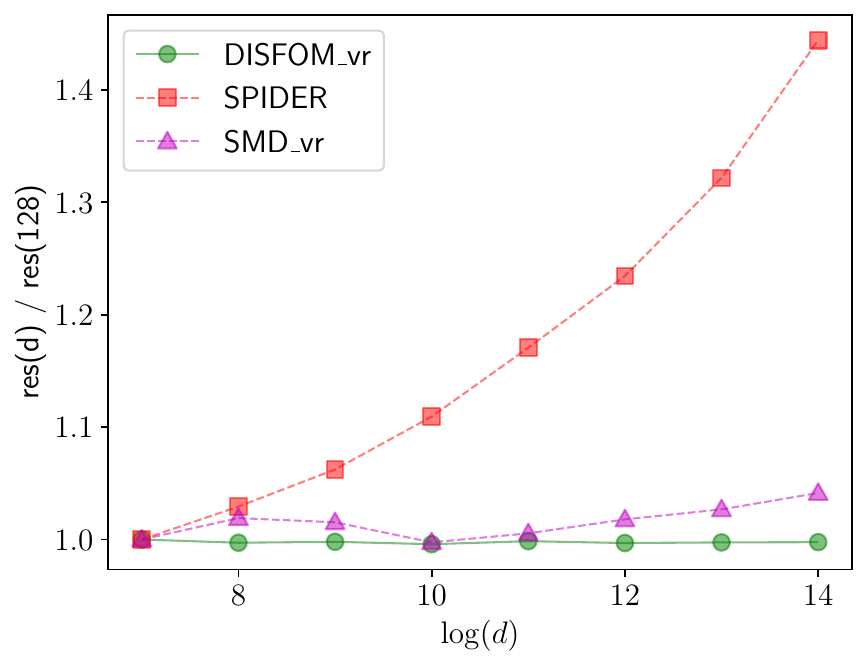}
    \caption{Averaged relative residual growth for minibatch (left) and variance reduction (right)}\label{fig:svrg}
\end{figure}

From Figure~\ref{fig:minibatch} and \ref{fig:svrg} we can see that the performance of vanilla SGD and SPIDER would deteriorate rapidly with the increase of $d$. On the contrary, the performance of DISFOMs and SMDs are comparable and less sensitive to $d$, in terms of both the averaged relative function value gap $f- f^*$ and the averaged relative residual defined in \eqref{def: resfunc}. 
To our best knowledge, there is not a proper theory for the dimension-insensitive property of SMD$\_$vr in literature, and the numerical experiments imply this empirically. 

\section{Conclusion}\label{sec: conclusion}

In this study, we introduce dimension-insensitive stochastic first-order methods (DISFOMs) as a solution for nonconvex optimization with an expected-valued objective function. Our algorithms are designed to accommodate non-Euclidean and non-smooth distance functions as proximal terms. When dealing with a stochastic nonconvex optimization problem, the sample complexity of stochastic first-order methods can become linearly dependent on the problem dimension. This causes trouble for solving large-scale problems. We demonstrate that DISFOM achieves a sample complexity of $O(L \Delta_f (\log d) \sigma_\infty^2 / \epsilon^4)$ in order to obtain an $\epsilon$-stationary point via minibatch sampling. Furthermore, we show that DISFOM with variance reduction can enhance this bound to $O(L \Delta_f (\log d)^2 \sigma_\infty / \epsilon^{3})$.  We present two options for the non-smooth distance functions, both of which allow for closed-form solutions to the subproblems in cases of unconstrained, box constraints and $\ell_1$-ball + box constraints. Preliminary numerical experiments demonstrate efficiency in solving the subproblems and illustrate the dimension-insensitive property of our proposed methods. Future directions include generalization to other non-Euclidean settings and consideration of non-smooth objective functions. 

\medskip
\noindent{\bf Acknowledgements} The research of Yue Xie is funded by Hong Kong RGC General Research Fund (Project number:  17300824),  and Guangdong Province Fundamental and Applied Fundamental Research Regional Joint Fund (Project number: 2022B1515130009).  The authors would like to thank Prof. George Lan, Prof. Anatoli Juditsky and Prof. Shuzhong Zhang for their valuable suggestions during writing this manuscript.

\section*{Statements and Declarations}
{\bf Competing Interests} The authors have no competing interests to declare that are relevant to the content of this article. \\
{\bf Availability of data and materials} The datasets generated during the current study are available in the GitHub repository, \\
\href{https://github.com/bjw010529/DISFOM.git}{https://github.com/bjw010529/DISFOM.git}.

		
\bibliographystyle{spmpsci}      
\bibliography{ref} 

	

\appendix

\section{Proofs}\label{app: proof}

\begin{lemma}\label{lm: subG}
    Suppose that $ \bar w \in \bR^d $ and $\bar w_i$ conforms sub-Gaussian distribution with norm $\sigma$. Then $\bE[ \| \bar w \|^2_\infty ] \le 6 \log (2d) \sigma^2$.
\end{lemma}
\begin{proof}
    According to the property of sub-Gaussian, we have that $\bE[ \exp(\bar w_i^2/(6\sigma^2))] \le 2$. Then, for any $t > 0$,
    \begin{align*}
       & \exp(t \bE[ \| \bar w \|_\infty^2]) \\
       & \le \bE[\exp( t \max_{1 \le j \le d} (\bar w_j)^2 )] \\
       & = \bE[ \max_{1 \le j \le d} \exp( t (\bar w_j)^2 )] \\
       & \le \bE[ \sum_{j = 1}^d \exp( t (\bar w_j)^2 )] \\
       & = \sum_{j = 1}^d \bE[ \exp( t (\bar w_j)^2 )],
    \end{align*}
    where the first inequality is from Jensen's inequality. If we let $t \textcolor{black}{\triangleq} 1/(6\sigma^2)$, then
    \begin{align*}
        & \exp\left( \frac{1}{6\sigma^2} \bE[ \| \bar w \|_\infty^2 ] \right) \le 2d \implies \bE[ \| \bar w \|_\infty^2 ] \le 6\log(2d) \sigma^2,
    \end{align*}
    and we proved the result.
    \qed
\end{proof}

Proof of Lemma~\ref{lm: l1subdiff}.
\begin{proof}
    Let $\conv(S)$ denote the convex hull of a set $S$. Note that $\| \x \|_1^2 = \max_{\alpha \in [-1,1]^d} (\alpha^T \x)^2 $. Then by Danskin's Theorem,
    {\small
    \begin{align*}
        & \partial \| \x \|_1^2 \\
        & = \conv\left\{ 2(\alpha^T \x)\alpha : \begin{cases}
            \alpha_i = 1 & \mbox{if } \x_i > 0, \\
            \alpha_i = -1 & \mbox{if } \x_i < 0, \\
            \alpha_i \in [-1,1] & \mbox{if } \x_i = 0,
        \end{cases} \mbox{ or }
        \begin{cases}
            \alpha_i = -1 & \mbox{if } \x_i > 0, \\
            \alpha_i = 1 & \mbox{if } \x_i < 0, \\
            \alpha_i \in [-1,1] & \mbox{if } \x_i = 0,
        \end{cases} 
        \qquad \right\}\\
        & = \conv\left\{ 2\| \x \|_1 \alpha : \begin{cases}
            \alpha_i = 1 & \mbox{if } \x_i > 0, \\
            \alpha_i = -1 & \mbox{if } \x_i < 0, \\
            \alpha_i \in [-1,1] & \mbox{if } \x_i = 0,
        \end{cases}  
        \qquad \right\} \\
        & = 2\| \x \|_1 \partial \| \x \|_1.
    \end{align*}}
    Therefore the result follows. \qed
\end{proof}

Proof of Theorem~\ref{thm: case1svrg}.
\begin{proof}
Note that similar to deriving \eqref{errbd3}, \eqref{errbd1}\eqref{setting1}\eqref{subgbd}\eqref{ineq: phi} and the parameter settings (let $t_k \equiv t$) lead to 
\begin{align}
\notag
& \sum_{k=1}^K \frac{1}{\eta} \bE [ \| \Delta \xbf^{k+1} \|^2 ] +  \sum_{k=1}^K \left( \frac{\hat \rho}{2 \eta } - \frac{L}{2} -  \frac{1}{2 t}  \right) \bE[ \| \Delta \xbf^{k+1} \|_1^2 ] \\
\notag
& \le \Delta_f +  \sum_{k=1}^K \frac{t}{2} \bE[ \| g^k - G^k \|_\infty^2 ] \\
\notag
& \implies \left( \frac{\hrho}{4\eta} - \frac{L}{2} \right) \bE [ \| \Delta \xbf^{Y+1} \|_1^2 ] \\\notag
& \le \frac{\Delta_f}{K} +  \frac{t}{2}  \bE[ \| g^Y - G^Y \|_\infty^2 ] \\\notag
& \overset{\eqref{redvar-SPIDER2}}{\le} \frac{\Delta_f}{K} +  \frac{t}{2}  \left( \frac{4 \Omega^2 L^2 q_0 }{m} \bE[ \| \x^{Y+1} - \x^Y \|_1^2 ] + \frac{\Omega e^2(\log d)\sigma_\infty^2}{m_1} \right) \\
& \implies \left( \frac{\hrho}{4\eta} - \frac{L}{2}  - \frac{2t \Omega^2 L^2 q_0}{m} \right) \bE [ \| \Delta \xbf^{Y+1} \|_1^2 ]
\label{errbd7}
\le \frac{\Delta_f}{K} + \frac{t \Omega e^2 (\log d)\sigmainf^2}{2 m_1}.
\end{align}
By \eqref{redvar-SPIDER2}\eqref{errbd2}\eqref{subgbd}\eqref{errbd7},
\begin{align}
\notag
    & \bE[ {\rm dist}_{\| \cdot \|_\infty}(0, \partial(f+\delta_X)(\x^{Y+1})) ] \\
\notag
& \le \left(L + \frac{1}{\eta} \right) \bE[ \| \Delta \x^{Y+1} \|_1 ] +  \bE[ \| g^Y - G^Y \|_\infty ] + \frac{1}{\eta} \bE[ \| \xi^{Y+1} \|_\infty ] \\ \notag
& \le \left(L + \frac{1}{\eta} \right) \sqrt{\bE[ \| \Delta \x^{Y+1} \|_1^2 ]} +  \sqrt{\bE[ \| g^Y - G^Y \|_\infty^2 ]} + \frac{1}{\eta} \bE[ \| \xi^{Y+1} \|_\infty ] \\ \notag
& \le \left(L + \frac{1}{\eta} \right) \sqrt{\bE[ \| \Delta \x^{Y+1} \|_1^2 ]} + \sqrt{ \frac{4 L^2 \Omega^2 q_0}{m} \bE[ \| \Delta \x^{Y+1} \|_1^2 ] + \frac{\Omega e^2 (\log d)\sigma_\infty^2}{m_1} } + \frac{ \hrho}{\eta} \bE[ \| \Delta \x^{Y+1} \|_1 ] \\ \notag
& \le \left(L + \frac{\hrho + 1}{\eta} + \sqrt{\frac{4L^2\Omega^2 q_0}{m}} \right) \sqrt{\bE[ \| \Delta \x^{Y+1} \|_1^2 ]} +  \sqrt{ \frac{\Omega e^2(\log d)\sigma_\infty^2}{m_1} } \\ \notag
    & \le \left( L + \frac{\hrho + 1}{\eta} + \sqrt{\frac{4\Omega^2 L^2q_0}{m}} \right) \sqrt{\frac{\frac{\Delta_f}{K} + \frac{ t \Omega e^2 (\log d)\sigmainf^2}{2m_1} }{\frac{\hrho}{4\eta} - \frac{L}{2} - \frac{2 t\Omega^2L^2q_0}{m} }} + \sqrt{ \frac{\Omega e^2 (\log d)\sigma_\infty^2}{m_1} }.
\end{align}
\qed
\end{proof}
Proof of Theorem~\ref{thm: case2svrg}.
\begin{proof}
By plugging in the parameter setting in this theorem and \eqref{errbd1}\eqref{setting2}\eqref{subgbd2}\eqref{ineq: nn3}, we have
\begin{align}
\notag
& \sum_{k=1}^K \frac{1}{\eta} \mathbb{E}[\|\Delta \x^{k+1}\|^2] - \sum_{k=1}^K \left(\frac{L}{2} + \frac{1}{2t}\right) \mathbb{E}[\|\Delta \x^{k+1}\|_1^2] + \sum_{k=1}^K \mathbb{E}\left[\frac{1}{\eta} \hat{\rho}_k \|\Delta \x^{k+1}\|_1\right] \\
\notag
&  \le{} \sum_{k=1}^K \frac{1}{\eta} \mathbb{E}[\|\Delta \x^{k+1}\|^2] - \sum_{k=1}^K \left(\frac{L}{2} + \frac{1}{2t_k}\right) \mathbb{E}[\|\Delta \x^{k+1}\|_1^2] + \sum_{k=1}^K \frac{1}{\eta} \mathbb{E}\left[(\xi^{k+1})^T \Delta \x^{k+1}\right] \\
\notag
&  \le{} \Delta_f + \sum_{k=1}^K \frac{t}{2} \mathbb{E}[\|g^k - G^k\|_\infty^2].
\end{align}
Therefore,
\begin{align}
\notag
&  \frac{\psi}{\eta} \bE[ \hrho_Y ] \\
\notag
& \le \frac{\Delta_f}{K} +  \frac{t}{2}  \bE[ \| g^Y - G^Y \|_\infty^2 ] + \left( \frac{L}{2} + \frac{1}{2 t} \right) \bE [ \| \Delta \xbf^{Y+1} \|_1^2 ] \\
\notag
& \overset{\eqref{redvar-SPIDER2}}{\le} \frac{\Delta_f}{K} + \frac{t}{2} \left( \frac{4 \Omega^2 L^2 q_0}{m} \bE[ \| \Delta \x^{Y+1} \|_1^2 ] + \frac{\Omega e^2 (\log d)\sigma_\infty^2}{m_1} \right) + \left( \frac{L}{2} + \frac{1}{2 t} \right) \bE [ \| \Delta \xbf^{Y+1} \|_1^2 ] \\
& \le \frac{\Delta_f}{K} +  \frac{ t \Omega e^2 (\log d)\sigma_\infty^2}{2 m_1} 
\label{errbd9}
+ \left( \frac{2 t  \Omega^2 L^2 q_0}{m} + \frac{L}{2} + \frac{1}{2 t} \right) \psi^2.
\end{align}
Then by \eqref{redvar-SPIDER2}\eqref{errbd2},
\begin{align}
\notag
    & \bE[ {\rm dist}_{\| \cdot \|_\infty}(0, \partial(f+\delta_X)(\x^{Y+1})) ] \\
\notag
& \le \left(L + \frac{1}{\eta} \right) \sqrt{\bE[ \| \Delta \x^{Y+1} \|_1^2 ]} +  \sqrt{\bE[ \| g^Y - G^Y \|_\infty^2 ]} + \frac{1}{\eta} \bE[ \| \xi^{Y+1} \|_\infty ] \\ \notag
& \le \left(L + \frac{1}{\eta} \right) \sqrt{\bE[ \| \Delta \x^{Y+1} \|_1^2 ]} + \sqrt{ \frac{4 \Omega^2 L^2q_0}{m} \bE[ \| \x^{Y+1} - \x^Y \|_1^2 ] + \frac{\Omega e^2 (\log d)\sigma_\infty^2}{m_1} } + \frac{1}{\eta} \bE[ \| \xi^{Y+1} \|_\infty ] \\ \notag
& \overset{\eqref{subgbd2}}{\le} \left(L + \frac{1}{\eta} + \sqrt{\frac{4 \Omega^2 L^2 q_0}{m}} \right) \sqrt{\bE[ \| \Delta \x^{Y+1} \|_1^2 ]} +  \sqrt{ \frac{\Omega e^2 (\log d)\sigma_\infty^2}{m_1} } + \frac{1}{\eta} \bE[  \hrho _Y ]
\\ \notag
    & \overset{\eqref{errbd9}}{\le} \left( L + \frac{1}{\eta} + \sqrt{\frac{4\Omega^2 L^2 q_0}{m}} \right) \psi + \sqrt{ \frac{\Omega e^2 (\log d)\sigma_\infty^2}{m_1} }  \\ 
    \notag
    & + \frac{1}{\psi} \left( \frac{\Delta_f}{K} + \frac{ t \Omega e^2 (\log d)\sigmainf^2}{2 m_1} + \left(\frac{2 t\Omega^2 L^2 q_0}{m} + \frac{L}{2} + \frac{1}{2t} \right) \psi^2 \right) .
\end{align}\qed
\end{proof}

Proof of Lemma~\ref{lm: l1consprox}.
\begin{proof}
The corresponding KKT-condition is as follows\begin{align}
\nabla_{\z} L(\z, \lambda) &= (\z - \vbf) + \lambda \partial\lVert \z \rVert_1 = 0 \label{eq:1}\\
g(x) &=  \lVert \z \rVert_1 - \psi \leq 0 \\
\lambda &\geq 0 \\
\lambda g(x) &= 0
\end{align}
From \eqref{eq:1} we have \begin{align*}
\z_i - \vbf_i + \lambda = 0 & \quad \text{if} \quad  \z_i > 0 \\
\z_i - \vbf_i - \lambda = 0 & \quad \text{if} \quad  \z_i < 0 \\
-\lambda \leq \z_i - \vbf_i \leq \lambda & \quad \text{if} \quad  \z_i = 0 \\
\end{align*}
Thus we have \begin{align}
\z_i = \vbf_i - \lambda & \quad \text{if} \quad \vbf_i - \lambda > 0 \label{eq:x1} \\
\z_i = \vbf_i + \lambda & \quad \text{if} \quad \vbf_i + \lambda < 0 \label{eq:x2} \\
\z_i = 0 & \quad \text{if} \quad \lvert \vbf_i \rvert \leq \lambda \label{eq:x3}
\end{align}
Suppose $\lambda = 0$, then $\z = \vbf$. This is the solution if $ \| \vbf \|_1 \le \psi$. Otherwise suppose that $\lambda > 0$. Then we have \begin{align} 
\lVert \z \rVert_1 = \sum_{i}^{\lvert \vbf_i \rvert > \lambda} \lvert \z_i \rvert 
= \sum_{i}^{\vbf_i > \lambda} (\vbf_i - \lambda) + \sum_{i}^{\vbf_i < -\lambda}(-\vbf_i - \lambda) = 
\sum_{i}^{\lvert \vbf_i\rvert > \lambda} (\lvert \vbf_i \rvert - \lambda) = \psi
\end{align}
Suppose there exists $m$ such that $s_{m-1} < \psi \le s_m$, then we can deduce that $| \vbf_{i_m} | > \lambda  \ge | \vbf_{i_{m+1}} |$, therefore,
\begin{align}
    \psi = \sum_{t=1}^m \lvert \vbf_{i_t}\rvert - m\lambda \implies 
    \lambda = \dfrac1{m}\left(\sum_{t=1}^{m} \left\lvert \vbf_{i_t} \right\rvert - \psi \right). \label{eq:v}
\end{align}
Finally we could derive $\z$ from (\ref{eq:x1}), (\ref{eq:x2}), (\ref{eq:x3}), (\ref{eq:v}) as follows.
\begin{align}
    \z_{i_t} = \begin{cases}
        0 & \mbox{if }  m+1 \le t \le d \\
        \vbf_{i_t} - \frac{\sgn\left(\vbf_{i_t}\right)}{m}\left(\sum_{k=1}^{m} \left\lvert \vbf_{i_k} \right\rvert - \psi \right) & \mbox{if }  1 \le t \le m
    \end{cases}
\end{align} \qed
\end{proof}

\section{Discussion on Lipschitz smoothness constant}\label{app: Lip}
{\revi
We supplement several examples to illustrate the advantage of using Lipschitz smoothness condition \eqref{ineq: lip} instead of its counterpart in the Euclidean setting.  Let us denote the Lipschitz constant in \eqref{ineq: lip} as $L_{\infty,1}$ and the one in the Euclidean setting as $L_{2,2}$.  For some problem classes,  $L_{\infty,1}$ can be much less dimension-sensitive than $L_{2,2}$.  Suppose that $f$ is $C^2$ over a neighborhood of $X$.  Note that for any $x,y \in X$ ($X$ is closed convex), we have
\begin{align*}
\| \nabla f(x) - \nabla f(y) \|_2 & = \| \nabla^2 f(y + \theta(x - y)) (x-y) \|_2 \\
& \le \| \nabla^2 f(y + \theta(x - y))\|_2 \| (x-y) \|_2,  \\
& \le  \max_{z \in X} \| \nabla^2 f(z)  \|_2 \| x-y \|_2,
\end{align*}
where $\theta$ is a certain scalar in $[0,1]$.  Therefore,  a typical estimate for $L_{2,2}$ is $\max_{x \in X} \| \nabla^2 f(x)  \|_2$. Meanwhile,  for any $x,y \in X$, we have
\begin{align*}
\| \nabla f(x) - \nabla f(y) \|_\infty & = \| \nabla^2 f(y + \theta(x - y)) (x-y) \|_\infty \\
& = \max_i \left| \sum_{j=1}^d \nabla^2_{ij} f(y + \theta(x-y))(x_j-y_j) \right| \\
& \le \max_i \left| \max_j | \nabla^2_{ij} f(y + \theta(x-y))| \| x-y \|_1 \right| \\
& = \max_{i,j} \left| \nabla^2_{ij} f(y + \theta(x-y)) \right| \| x-y \|_1  \\
& \le  \max_{z \in X,i,j} \left| \nabla^2_{ij} f(z) \right|  \| x-y \|_1. 
\end{align*}
Then we will let $L_{\infty,1} = \max_{z \in X,i,j} \left| \nabla^2_{ij} f(z) \right| $.  

Let us consider two problem classes.  The first one is when $f$ is a quadratic function and suppose that its Hessian $H$ has randomly generated entries.  Specifically,  $H$ is symmetric and suppose that its entries in upper triangular region and on the diagonal conform i.i.d uniform distribution on $[-a,a]$.  Then by \cite[Theorem A]{bai1988necessary}, we have $\lambda_{max}(H /\sqrt{d}) \to 2 \sqrt{a/3}$ and $\lambda_{max}(-H /\sqrt{d}) \to 2 \sqrt{a/3}$ almost surely when $d \to \infty$. Therefore,  $\| H \|_2 = 2\sqrt{a d / 3 } + o(\sqrt{d})$ when $d\to \infty$ almost surely.  On the other hand,  $| \max_{i,j} H_{i,j}| \le a$ with probability $1$ for any $d$.  The former is significantly larger and dimension sensitive than the latter.

Another problem class is when $f(x) \triangleq h( a^T x - b)$, $ a \in \mathbb{R}^d$, $b \in \mathbb{R}$ where $h$ is a certain loss function,  e.g.,  square loss,  logistic loss,  smooth nonconvex surrogates for $\ell_0$ norm, etc.  Therefore,  $\nabla^2 f(x) = h''(a^T x - b) a a^T $.  Suppose that $\lambda \triangleq \sup_{z \in Z} |h''(z) |$ where $Z = \{ z \mid z = a^T x - b, x \in X \}$.  For e.g.,  if $h(z) = (1/2)z^2$, then $\lambda = 1$; if $h(z) = \log(1 + e^{-z})$ and $0 \in Z$, then $h''(z) = e^z/(1 + e^z)^2$ and $\lambda = 1/4$.  Therefore, we have $\max_{x \in X} \|  \nabla^2 f(x) \|_2 = \lambda  \| a a^T \|_2 = \lambda \| a \|_2^2 $, which could grow linear in $d$.  On the other hand,  $\max_{z \in X,i,j} \left| \nabla^2_{ij} f(z) \right| = \lambda \max_{i,j} | a_i a_j | = \lambda \| a \|_\infty^2$,  which is much less dimension-sensitive. 

For the problem classes discussed above,  the estimate of $L_{2,2}$ can be significantly larger than $L_{\infty,1}$ when $d$ is large.  
}

\section{Lower bounds}\label{app: lb}
We follow the settings in \cite{arjevani2023lower}. Suppose that for the rest of this section, $X = \bR^d$, and we denote the function class that satisfies Assumption~\ref{asp: flip},~\ref{ass: setw3hj} as $$
    \mathcal{H}_\infty(L, \Delta_f) \triangleq \{ f: f(0) - \inf f \le \Delta_f, \| \nabla f(\x) - \nabla f(\y) \|_\infty \le L \| \x - \y \|_1, \forall \x, \y \}. $$
    The algorithm accesses the information of $f \in \mathcal{H}_\infty(L, \Delta_f)$ via a stochastic first-order oracle $O$. Suppose that at each sampling round $i$, the algorithm considered queries a batch of size $K$, i.e., 
    $\x^{(i)} \textcolor{black}{\triangleq} \{ \x^{(i,1)},...,\x^{(i,K)} \}$, where $\x^{(i,k)} \in \bR^d$ and $k \in [K]$. For each batch query $\x^{(i)}$, the oracle $O$ performs an independent draw $z^{(i)} \sim P_z$ and responds with $$O_f(\x^{(i)}, z^{(i)}) \textcolor{black}{\triangleq} (O_f(\x^{(i,1)}, z^{(i)}),...,O_f(\x^{(i,K)}, z^{(i)})).$$ Variance reduction corresponds to $K = 2$. 
    $\mathcal{O}_{\infty}(K,\sigma_{\infty}^2)$ denotes the oracle class with unbiasedness and bounded variance in $\| \cdot \|_{\infty}$ norm (Assumption~\ref{ass: subG}). $\mathcal{O}_{\infty}(K,\sigma_{\infty}^2, L)$ denotes the oracle class that additionally satisfies mean-square smoothness (Assumption~\ref{ass: Lip}). An algorithm $A$ consists of a distribution $P_r$ over a measurable set $\mathcal{R}$ and a sequence of measurable mappings $\{ A^{(i)} \}_{i \in \mathbb{N}}$ such that $A^{(i)}$ takes in the first $i-1$ oracle responses and the random seed $r \in \mathcal{R}$ to produce the $i$th query. Let $\{ \x_{A[O_f]}^{(i)} \}_{i \in \mathbb{N}}$ denotes the random sequence of queries resulting from applying $A$ with $O$, defined recursively as
    \begin{align}\label{optpro}
        \x_{A[O_f]}^{(i)} \textcolor{black}{\triangleq} A^{(i)}\left(r, O_f\left(\x_{A[O_f]}^{(1)}, z^{(1)}\right),...,O_f\left(\x_{A[O_f]}^{(i-1)}, z^{(i-1)}\right)\right),
    \end{align}
    where $r \sim P_r$ is drawn a single time at the beginning of the optimization protocol. Use $\sA(K)$ to denote the class of algorithms that follow this protocol \eqref{optpro}. Then the \textit{distributional complexity} to find an $\epsilon$-stationary solution \eqref{def: epstat}\eqref{def: resfunc} is defined as
    \begin{align}\label{def: dist.comp.lip}
        m_{\infty,\epsilon} (K,\sigma_{\infty}^2,L,\Delta_f) & \triangleq \sup_{O \in \sO_{\infty}(K,\sigma_\infty^2)} \;  \sup_{P_f \in \sP(\sH_\infty(L,\Delta_f))} \; \inf_{A \in \sA(K)} \\
        \notag
        & \inf \{ T \mid \bE[ \| \nabla f(\x_{A[O_f]}^{(T)}) \|_\infty ] \le \epsilon \},
    \end{align}
    where $\sP$ denotes the class of all the distributions over $\sH_\infty(L,\Delta_f)$. In the mean-squared smoothness setting, this complexity is defined as 
    \begin{align}\label{def: dist.comp.ms}
        \bar{m}_{\infty,\epsilon} (K,\sigma_{\infty}^2,L,\Delta_f) & \triangleq \sup_{O \in \sO_{\infty}(K,\sigma_\infty^2,L)} \;  \sup_{P_f \in \sP(\sH_\infty(L,\Delta_f))} \; \inf_{A \in \sA(K)} \\
        \notag
        & \inf \{ T \mid \bE[ \| \nabla f(\x_{A[O_f]}^{(T)}) \|_\infty ] \le \epsilon \}.
    \end{align}
Sample complexity results in Corollary~\ref{corr: mb.comp.case1},~\ref{corr: mb.comp.case2} and Corollary~\ref{corr: vr.comp.case1},~\ref{corr: vr.comp.case2} can be regarded as upper bounds for \eqref{def: dist.comp.lip} and \eqref{def: dist.comp.ms} respectively. Now we describe the lower bounds. We restrict to the zero-respecting algorithm class, which is a fairly general class including many first-order algorithm designs.  
\begin{definition}
    A stochastic first-order algorithm $A$ is zero-respecting if for any oracle $O$ and any realization of $z^{(1)}, z^{(2)}, ...$, for all $t \ge 1$ and $k \in [K]$, 
    \begin{align*}
        \mbox{support}(\x_{A[O_f]}^{t,k}) \subseteq \bigcup_{i < t, k' \in [K]} \mbox{support}( g^{(i,k')} ),
    \end{align*}
    where $(f^{(t,1)},g^{(t,1)}),...,(f^{(t,K)},g^{(t,K)}) \textcolor{black}{\triangleq} O_f(\x_{A[O_f]}^{(t)}, z^{(t)})$ denote the oracle responses for round $t$. Let $\sA_{zr}(K)$ denote the class of zero-respecting algorithms.
\end{definition}
Following are lower bound results for the distributional complexity of $\sA_{zr}(K)$.
\begin{theorem}\label{thm: lb1}
    There exists numerical constants $c,c' > 0$ such that for all $L$, $\Delta_f$, $\sigma_\infty^2 > 0$ and $\epsilon \le c' \sqrt{L \Delta_f}$,
    \begin{align*}
        m^{zr}_{\infty,\epsilon}(K,\Delta_f,L,\sigma_\infty^2) \ge c \left( \frac{\Delta_f L \sigma_\infty^2}{\epsilon^4} + \frac{\Delta_f L}{ \epsilon^2 } \right).
    \end{align*}
    Constructions of dimension $d = \mathcal{O}\left( \frac{\Delta_f L}{\epsilon^2} \right)$ realize the lower bound.
\end{theorem}
\begin{theorem}\label{thm: lb2}
    There exist numerical constants $c$, $c' > 0$ such that for all $L$, $\Delta_f$, $\sigma_\infty^2 > 0$ and $\epsilon \le c' \sqrt{L \Delta_f}$, 
    \begin{align*}
        \bar{m}^{zr}_{\infty,\epsilon}(K,\Delta_f,L,\sigma_\infty^2) \ge c\left( \frac{\Delta_f L \sigma_\infty}{\epsilon^3} + \frac{\Delta_f L}{\epsilon^2} + \frac{\sigma_\infty^2}{\epsilon^2} \right).
    \end{align*}
    Construction of dimension $d = \mathcal{O}\left( 1 + \frac{\Delta_f L}{\sigma_\infty \epsilon} \right)$ realize the lower bound.
\end{theorem}
\begin{proof}
    Proofs of Theorem~\ref{thm: lb1} and Theorem~\ref{thm: lb2} follow closely the proofs of Theorem 1 and 2 in \cite{arjevani2023lower}. The same hard examples accommodate the $\| \cdot \|_1$,$\| \cdot \|_\infty$ setting. \qed
\end{proof}

\section{Algorithm}\label{app: alg}
\paragraph{Gradient Descent with Backtracking}\quad
Here we provide the details of the gradient descent with backtracking (Algorithm \ref{alg: GD backtracking}) used for computing $f^*$ in the numerical experiment. 

\begin{algorithm} 
\caption{Projected gradient method with backtracking}\label{alg: GD backtracking}
\begin{description}
\item[{\bf Initialization Step}] Set $c_1 \textcolor{black}{\triangleq} \dfrac{1}{4}$, $\beta \textcolor{black}{\triangleq} \dfrac{1}{2}$, and $\epsilon \textcolor{black}{\triangleq} 1e-10$. Initialize $\xbf^0 \textcolor{black}{\triangleq} \bm{0}$.
\item[{\bf Main Step}] Set $\alpha_k \equiv \alpha \textcolor{black}{\triangleq} 1$;\\
{\bf While} {$f(P_X (\xbf^k - \alpha_k \nabla f(\xbf^k))) > f(\xbf^{k}) + c_1 (\nabla f(\xbf^{k}))^T (P_X (\x^k - \alpha_k \nabla f(\x^k)) - \x^k) $}\\
$\alpha_k = \beta\alpha_k$
\item[{\bf Stopping}] Stop if $\lVert \xbf^{k+1} - \xbf^{k} \rVert_1 \le \epsilon$ and output $\xbf^{k}$; otherwise replace $k$ by $k+1$ and repeat the Main Step.
\end{description}
\end{algorithm}

\end{document}